\documentclass[11pt]{article}\textwidth 160mm\textheight 235mm
\usepackage{amsfonts}
\usepackage{amssymb}
\usepackage{graphicx}
\usepackage{amsmath}
\usepackage{tikz}
\usepackage{cancel}
\usepackage{todonotes}
\usetikzlibrary{matrix}
\usetikzlibrary{arrows,shapes}

\def\ve{\varepsilon}
\def\epsilon{\varepsilon}

\def\hat{\widehat}
\def\tilde{\widetilde}
\def\emp{\emptyset}

\def\dist{{\rm dist}}
\def\dom{{\rm dom}\,}
\def\rge{{\rm rge}\,}

\def\rge{{\rm rge\,}}

\def\N{{\cal N}}

\def\G{{\Gamma}}
\def\O{{\cal O}}

\def\sub{\partial}
\def\B{\mathbb B}
\def\Rm{{\mathbb R}^{m+1}}

\def\ox{\overline{x}}
\def\oy{\overline{y}}
\def\oz{\overline{z}}

\def\disp{\displaystyle}

\def\Limsup{\mathop{{\rm Lim}\,{\rm sup}}}

\def\tto{\;{\lower 1pt \hbox{$\rightarrow$}}\kern-10pt
\hbox{\raise 2pt\hbox{$\rightarrow$}}\;}
\def\Hat{\widehat}

\def\Bar{\overline}
\def\ra{\rangle}
\def\la{\langle}
\def\ve{\varepsilon}
\def\B{I\!\!B}
\def\h{\hfill\Box}
\def\R{\mathbb{R}}
\def\N{I\!\!N}
\def\ox{\bar{x}}
\def\oy{\bar{y}}
\def\oz{\bar{z}}

\def\ou{\bar{u}}
\def\ot{\bar{t}}
\def\op{\bar{p}}

\def\ri{\mbox{\rm ri}}
\def\int{\mbox{\rm int}}
\def\gph{\mbox{\rm gph}\,}

\def\dom{\mbox{\rm dom}\,}
\def\bd{\mbox{\rm bd}}
\def\ker{\mbox{\rm ker}\,}

\def\h{\hfill\triangle}

\def\O{\Omega}
\def\ph{\varphi}
\def\emp{\emptyset}
\def\st{\stackrel}

\def\lm{\lambda}
\def\olm{\bar\lambda}
\def\gg{\gamma}
\def\dd{\delta}
\def\al{\alpha}
\def\Th{\Theta}
\def\N{I\!\!N}

\def\sce{\setcounter{equation}{0}}
\def\Q{{\cal Q}}

\def\H{{\cal H}}
\begin{document}
\vspace*{0.5in}
\begin{center}
{\bf SECOND-ORDER VARIATIONAL ANALYSIS\\IN SECOND-ORDER CONE PROGRAMMING}\\[2ex]
NGUYEN T. V. HANG\footnote{Department of Mathematics, Wayne State University, Detroit, MI 48202, USA and Institute of Mathematics, Vietnam Academy of Science and Technology, Hanoi 10307, Vietnam (hangnguyen@wayne.edu). Research of this author was partly supported by the National Science Foundation under grant DMS-1512846 and by the Air Force Office of Scientific Research under grant \#15RT0462.}, BORIS S. MORDUKHOVICH\footnote{Department of Mathematics, Wayne State University, Detroit, MI 48202, USA and RUDN University, Moscow 117198, Russia (boris@math.wayne.edu). Research of this author was partly supported by the National Science Foundation under grant DMS-1512846, by the Air Force Office of Scientific Research under grant \#15RT0462, and by the Ministry of Education and Science of the Russian Federation (Agreement number 02.a03.21.0008 of 24 June 2016).} and M. EBRAHIM SARABI\footnote{Department of Mathematics, Miami University, Oxford, OH 45065, USA (sarabim@miamioh.edu).}
\end{center}
\vspace*{0.05in}

\small{\bf Abstract.} The paper conducts a second-order variational analysis for an important class of nonpolyhedral conic programs generated by the so-called second-order/Lorentz/ice-cream cone $\Q$. From one hand, we prove that the indicator function of $\Q$ is always twice epi-differentiable and apply this result to characterizing the uniqueness of Lagrange multipliers at stationary points together with an error bound estimate in the general second-order cone setting involving ${\cal C}^2$-smooth data. On the other hand, we precisely calculate the graphical derivative of the normal cone mapping to $\Q$ under the weakest metric subregularity constraint qualification and then give an application of the latter result to a complete characterization of isolated calmness for perturbed variational systems associated with second-order cone programs. The obtained results seem to be the first in the literature in these directions for nonpolyhedral problems without imposing any nondegeneracy assumptions.\\[1ex]

\newtheorem{Theorem}{Theorem}[section]
\newtheorem{Proposition}[Theorem]{Proposition}
\newtheorem{Remark}[Theorem]{Remark}
\newtheorem{Lemma}[Theorem]{Lemma}
\newtheorem{Corollary}[Theorem]{Corollary}
\newtheorem{Definition}[Theorem]{Definition}
\newtheorem{Example}[Theorem]{Example}
\renewcommand{\theequation}{{\thesection}.\arabic{equation}}
\renewcommand{\thefootnote}{\fnsymbol{footnote}}

\normalsize

\section{Introduction}\sce

This paper is devoted to the study of {\em second-order cone programs} (SOCPs), which are optimization problems with constraint sets given by
\begin{equation}\label{CS}
\Gamma:=\big\{x\in\mathbb{R}^n\,\big|\,\Phi(x)\in\Q\,\big\},
\end{equation}
where $\Phi:\R^n\to\R^{m+1}$ is twice differentiable $({\cal C}^2$-smooth) around the reference points, and where ${\cal Q}$ is the {\em second-order/Lorentz/ice-cream cone} defined by
\begin{equation}\label{Qm}
\Q:=\big\{s=(s_0,s_r)\in\R\times\R^{m}\,\big|\,\Vert s_r\Vert\le s_0\big\}.
\end{equation}
Problems of this type are challenging mathematically while being important for various applications; see, e.g., \cite{ag,BoR05,bs,mos,or,os} and the bibliographies therein. A remarkable feature of SOCPs, which significantly distinguishes them from nonlinear programs (NLPs) and the like, is the {\em nonpolyhedrality} of the underlying Lorentz cone $\Q$ in \eqref{CS}.

The main intention of this paper is to conduct a second-order analysis for SOCPs by using appropriate tools of second-order generalized differentiation. The following two topics are of our particular interest here: (1) to prove the {\em twice epi-differentiability} of the indicator function of $\Q$ with deriving an explicit formula for the calculation of the second epi-derivative; (2) to establish a precise formula for {\em calculating the graphical derivative} of the {\em normal cone mapping} generated by the constraint set $\G$ in \eqref{CS} {\em without} imposing any {\em nondegeneracy} condition. To the best of our knowledge, the results obtained in both directions are the {\em first ones} in the literature for nonpolyhedral systems. They have strong potentials for applications to SOCPs and related problems. Among those presented in this paper we mention characterizations of the {\em uniqueness of Lagrange multipliers} together with an {\em error bound} estimate in SOCPs and also of the {\em isolated calmness} property for solution maps of perturbed variational systems associated with SOCPs.

Although for brevity and simplicity we consider here the case of one Lorentz cone \eqref{Qm} in \eqref{CS}, the developed approaches and results can be easily extended to the {\em product} setting
\begin{equation*}
{\cal Q}\colon=\prod_{j=1}^{J}{\cal Q}_{m_j+1}\subset\R^l\;\;\mbox{ with }\;l:=\sum_{j=1}^J(m_j+1),
\end{equation*}
where each cone ${\cal Q}_{m_j+1}$ is defined as in \eqref{Qm}; see Remark~\ref{product} for more details.\vspace*{0.03in}

The rest of the paper is organized as follows. In Section~2 we briefly overview the basic variational notions and constructions widely used in the sequel.

Section~3 concerns {\em second-order epi-differentiability} (in the sense of Rockafellar \cite{roc}) of the indicator function $\dd_\Q$ of the Lorentz cone \eqref{Qm}, some of its consequences, and related properties. The main result here not only justifies the twice epi-differentiability of $\dd_\Q$, but also establishes a precise formula for calculating the second-order epi-derivative of this function in terms of the given data of the Lorentz cone $\Q$ without any additional assumptions.

In Section~4 we study second-order properties of the SOCP constraint system \eqref{CS} by using the twice epi-differentiability of $\dd_\Q$ and the {\em metric subregularity constraint qualification} (MSCQ) for \eqref{CS}, which seems to be the weakest constraint qualification that has been investigated and employed recently in the (polyhedral) NLP framework; see \cite{go16,gm15,ch}. Among the most important results obtained in this section we mention the following: (i) a constructive description of generalized normals to the critical cone at the point in question under MSCQ, and (ii) a characterization of the uniqueness of Lagrange multipliers {\em together} with an appropriate error bound estimate (automatic in the polyhedral case) at stationary points via a new constraint qualification in conic programming, which happens to be in the case of \eqref{CS} a {\em dual} form of the {\em strict Robinson constraint qualification} (SRCQ) from \cite{bs,dsz}. We also present here novel {\em approximate duality} relationships for a linear conic optimization problem associated with the second-order cone $\Q$ that play a significant role in establishing the main result of the paper.

Section~5 derives a new formula allowing us to precisely calculate the {\em graphical derivative} of the normal cone mapping generated by \eqref{CS}, merely under the validity of MSCQ. The obtained major result is the first in the literature for nonpolyhedral constraint systems without imposing nondegeneracy. As discussed below, its proof is significantly different from the recent ones given in \cite{ch,gm15,go16} for polyhedral systems, even in the latter case. It is also largely different from the approaches developed in \cite{go17,mor15,mor} for conic programs under nondegeneracy assumptions.

In Section~6 we present a nontrivial example of a two-dimensional constraint system \eqref{CS} with the tree-dimensional Lorentz cone $\Q$ illustrating applications of the graphical derivative formula from Section~5. In this example the MSCQ condition holds at any feasible point of \eqref{CS} while the nondegeneracy and metric regularity/Robinson constraint qualification fail therein. We also apply the obtained graphical derivative formula to deriving a complete {\em characterization} of the {\em isolated calmness} property for solution maps to canonically perturbed variational systems associated with SOCP and give a numerical example.

The concluding Section~7 contains some discussions on further developments and applications of the approach and results of this paper in conic programming.\vspace*{0.03in}

Our notation and terminology are standard in variational analysis, conic programming, and generalized differentiation; see, e.g., \cite{bs,m06,rw}. Recall that, given a nonempty set $\O\subset\R^n$, the symbol $x\st{\O}{\to}\ox$ indicates that $x\to\ox$ with $x\in\O$. We often write an element $x\in\Rm$ in the second-order cone $\Q$ as $x=(x_0,x_r)$ with $x_0\in \R$ and $x_r\in \R^m$. Taking into account this decomposition of $x\in\Q$, denote $\hat{x}:=(-x_0,x_r)$. Finally,
$\B$ stands for the closed unit ball in the space in question while $\B_\gg(x):=x+\gg\B$ is the closed ball centered at $x$ with radius $\gg>0$.\vspace*{-0.1in}

\section{Preliminaries from Variational Analysis}\sce\vspace*{-0.05in}

In this section we first recall, following mainly the books \cite{m06,rw}, some basic notions from variational analysis and generalized differentiation widely used in the paper and then formulate the needed description of twice epi-differentiability for extended-real-valued functions.

Given a nonempty set $\O\subset\R^n$ locally closed around $\ox\in\O$, the (Bouligand-Severi) {\em tangent/contingent cone} $T_\O(\ox)$ to $\O$ at $\ox\in\O$ is defined by
\begin{eqnarray}\label{2.5}
T_\O(\ox):=\big\{w\in\R^n\big|\;\exists\,t_k{\downarrow}0,\;\;w_k\to w\;\;\mbox{ as }\;k\to\infty\;\;\mbox{with}\;\;\ox+t_kw_k\in \O\big\}
\end{eqnarray}
while the (Mordukhovich) {\em basic/limiting normal cone} to $\O$ at this point is given by
\begin{eqnarray}\label{nc}
N_\O(\ox)=\disp\Limsup_{x\to\ox}\Big[\mbox{cone}\big(x-\Pi_\O(x)\big)\Big],
\end{eqnarray}
where $\Pi_\O\colon\R^n\tto\R^n$ stands for the {\em Euclidean projector}. Despite the nonconvexity of the limiting normal cone \eqref{nc}, it enjoys---together with the associated subdifferential and coderivative constructions for extended-real-valued functions and set-valued mappings/multifunctions, respectively,---comprehensive calculus rules based on variational/extremal principles of variational analysis; see \cite{m06,rw} for more details. Note that a rich calculus is not available for the contingent cone \eqref{2.5}. If the set $\O$ is convex, then constructions \eqref{2.5} and \eqref{nc}  reduce, respectively, to the classical tangent and normal cones of convex analysis.

Considering next a set-valued mapping $F\colon\R^n\tto\R^m$ with its domain and graph given by
$$
\dom F:=\big\{x\in\R^n\big|\;F(x)\ne\emp\big\}\quad\mbox{and}\quad\gph F:=\big\{(x,y)\in\R^n\times\R^m\big|\;x\in F(x)\big\},
$$
we define the following generalized differential notions for $F$ induced by by the above tangent and normal cones to it graph. Given $(\ox,\ot)\in\gph F$, the {\em graphical derivative} of $F$ at $(\ox,\oy)$ is
\begin{equation}\label{gder}
DF(\ox,\oy)(u):=\big\{v\in\R^m\big|\;(u,v)\in T_{\small{\gph F}}(\ox,\oy)\big\},\quad u\in\R^n,
\end{equation}
while the {\em limiting coderivative} to $F$ at $(\ox,\oy)$ is defined by
\begin{equation}\label{2.8}
D^*F(\ox,\oy)(v):=\big\{u\in\R^n\big|\;(u,-v)\in N_{\small{\gph F}}(\ox,\oy)\big\},\quad v\in\R^m.
\end{equation}

Recall that a set-valued mapping $F\colon\R^n\tto\R^m$ is {\em metrically regular} around $(\ox,\oy)\in\gph F$ if there is $\ell\ge 0$ such that we have the distance estimate
\begin{equation}\label{metreq}
\dist\big(x;F^{-1}(y)\big)\le\ell\,\dist\big(y;F(x)\big)\quad\quad\mbox{for all}\quad (x,y)\;\; \mbox{close to} \;\;(\ox,\oy).
\end{equation}
If $y=\oy$ in \eqref{metreq}, the mapping $F$ is called to be {\em metrically subregular} at $(\ox,\oy)$.

Given $\O\subset\R^n$ and its indicator function $\delta_\O(x)$ equal to $0$ for $x\in\O$ and $\infty$ for $x\notin\O$, consider the parametric ($t>0)$ family
second-order difference quotients
\begin{equation}\label{lk01}
\Delta_t^2\delta_\O(\ox|\oy)(v):=\dfrac{\delta_\O(\ox+tv)-\delta_\O(\ox)-t\langle \oy,\,v\rangle}{\frac{1}{2}t^2}\quad\quad\mbox{with}\;\;v\in \R^{n},
\end{equation}
and say that $\delta_\O$ is {\em twice epi-differentiable} at $\bar x\in\O$ for $\oy\in\R^n$ in the sense of Rockafellar with its second-order epi-derivative $d^2\delta_\O(\bar x|\oy)\colon\R^n\to\overline{\mathbb{R}}$ if the second-order difference quotients $\Delta_t^2\delta_\O(\ox|\oy)$ epi-converges to $d^2\delta_\O(\bar x|\oy)$ as $t\downarrow 0$. The latter means by \cite[Proposition~7.2]{rw} that $\Delta_t^2\delta_\O(\bar x|\oy)$
for every sequence $t_k\downarrow 0$ and every $v\in\R^n$ we have
$$
\begin{array}{lll}
d^2\delta_\O(\bar x|\oy)(v)&\le&\disp\liminf_{k\to\infty}\Delta_{t_k}^2\delta_\O(\bar x|\oy)(v_k)\quad\mbox{for every sequence}\;\;v_k\to v,\\
d^2\delta_\O(\bar x|\oy)(v)&\ge&\disp\limsup_{k\to\infty}\Delta_{t_k}^2\delta_\O(\bar x|\oy)(v_k)\quad\mbox{for some sequence}\;\;v_k\to v.
\end{array}
$$\vspace*{-0.25in}

\section{Twice Epi-Differentiability of the Indicator Function of $\Q$}\sce\label{sec3}\vspace*{-0.05in}

We begin our second-order analysis with the study of twice epi-differentiability of the indicator function $\dd_\Q$ of the second-order cone \eqref{Qm}. The notions of first- and second-order epi-differentiability for extended-real-valued functions was introduced by Rockafellar in \cite{roc}, where he proved the twice epi-differentiability of convex piecewise linear-quadratic functions in finite dimensions. This result was extended in \cite[Theorem~14.14]{rw} to the class of fully amenable functions based on their polyhedral structure. Furthermore, Levy \cite{Lev99} established the twice epi-differentiability of the indicator function of convex polyhedric sets in Banach spaces; the latter notion reduces to the standard polyhedrality for sets in finite dimensions.

The following theorem justifies the twice epi-differentiability of the indicator function $\dd_\Q$ of the second-order cone \eqref{Qm} and calculates its second-order epi-derivative via the given data of $\Q$ without any additional assumptions. It seems to be the first result in the literature on second-order epi-differentiability in nonpolyhedral settings.\vspace*{-0.1in}

\begin{Theorem}{\bf(second-order epi-derivative of the indicator function of $\Q$).}\label{ted} Given any $\ox\in\Q$, the indicator function $\dd_\Q$ is twice epi-differentiable at $\ox$ for every $\oy\in N_\Q(\ox)$ and its second-order epi-derivative is calculated by
\begin{equation}\label{sed}
d^2\delta_\Q(\ox|\oy)(v)=\left\{\begin{array}{ll}
0\quad&\mbox{if }\;\ox\in[\int(\Q)\cup\{0\}],\;v\in\Bar{\cal K},\\
\disp\frac{\|\oy\|}{\|\ox\|}(\Vert v_r\Vert^2-v_0^2)\quad
&\mbox{if }\;\ox\in\bd(\Q)\setminus\{0\},\;v\in{\Bar{\cal K}},\\
\infty&\mbox{if }\;v\notin{\Bar{\cal K}},
\end{array}
\right.
\end{equation}
where ${\Bar{\cal K}}:=T_\Q(\ox)\cap \{\oy\}^\bot$ is the critical cone of $\Q$ at $\bar x$ for $\bar y$.
\end{Theorem}\vspace*{-0.1in}
{\bf Proof.} Fix $\ox\in\Q,\,\oy\in N_\Q(\ox)$, and $v\in\mathbb{R}^{m+1}$ and denote by $\Delta(\ox,\, \oy)(v)$  the right-hand side of \eqref{sed}.
To verify formula (\ref{sed}), we apply \cite[Proposition~7.2]{rw} that gives us the following description of the twice epi-differentiability of $\dd_\Q$ at $\ox$ for $\oy$:
\begin{itemize}
\item For every sequences $t_k\downarrow 0$ and $v_k\to v$ the second-order difference quotients \eqref{lk01} satisfy
\begin{equation}\label{sed1}
\liminf_{k\to\infty}\Delta_{t_k}^2\delta_\Q(\ox|\oy)(v_k)\ge\Delta(\ox|\oy)(v).
\end{equation}
\item For every sequence $t_k\downarrow 0$ there is some sequence $v_k\to v$ satisfying the inequality
\begin{equation}\label{sed2}
\limsup _{k\to\infty}\Delta_{t_k}^2\delta_\Q(\ox|\oy)(v_k)\le\Delta(\ox|\oy)(v).
\end{equation}
\end{itemize}
We split the proof into considering the three cases for $\ox\in\Q$ in representation \eqref{sed}.

{\bf Case~1:} $\ox\in\int(\Q)$. In this case we have $N_\Q(\ox)=\{0\}$ and hence $\oy=0$. Fix $v\in {\Bar{\cal K}}=\mathbb{R}^{m+1}$ and observe from (\ref{sed}) that $\Delta(\ox,0)(v)=0$. Picking an arbitrary sequence $v_k\to v$ as $k\to\infty$, we arrive at the equalities
\begin{equation*}
\Delta_{t_k}^2\delta_\Q(\ox|0)(v_k)=\dfrac{\delta_\Q(\ox+t_kv_k)-\delta_\Q(\ox)-t_k\cdot 0}{\frac{1}{2}t_k^2}=0
\end{equation*}
for all $k$ sufficiently large. This tells us that
\begin{equation*}
\lim_{k\to 0}\Delta_{t_k}^2\delta_\Q(\ox|0)(v_k)=0=\Delta(\ox,0)(v),
\end{equation*}
which confirms the validity of (\ref{sed1}) and (\ref{sed2}) and thus justifies formula (\ref{sed}) in this case.

{\bf Case~2:} $\ox=0$. In this case we clearly have $\oy\in N_\Q(\ox)=-\Q$. Pick $v\in\mathbb{R}^{m+1}$ and let $v_k\to v$ as $k\to\infty$. Using (\ref{lk01}) gives us the representations
\begin{equation}\label{eq16}
\Delta_{t_k}^2\delta_\Q(0|\oy)(v_k)=\dfrac{\delta_\Q(t_kv_k)-\delta_\Q(0)-t_k\langle \oy,v_k\rangle}{\frac{1}{2}t_k^2}=\begin{cases}
-\dfrac{\langle\oy,v_k\rangle}{\frac{1}{2}t_k^2}\ge 0\quad&\textrm{if }\;v_k\in \Q,\\
\infty&\textrm{if }\;v_k\notin \Q.
\end{cases}
\end{equation}
If $v\in{\Bar{\cal K}}$, we conclude from the above definition $\Delta(\ox,\oy)(v)$ that $\Delta(0,\oy)(v)=0$. Thus \eqref{sed1} comes directly from \eqref{eq16}, while (\ref{sed2}) can be justified  by choosing $v_k=v$ for any $k$. Pick now  $v\notin{\Bar{\cal K}}=\Q\cap\{\oy\}^\perp$ and observe that it amounts to saying that either $v\notin\Q$ or $\langle\oy,v\rangle<0$. It follows from the definition od $\Delta(\ox,\oy)(v)$ in this case that $\Delta(0,\oy)(v)=\infty$, and hence inequality (\ref{sed2}) clearly holds. To verify (\ref{sed1}), pick an arbitrary sequence $v_k\to v$. If $v\notin \Q$, then we can assume without loss of generality that $v_k\notin\Q$ for all $k$, which together with (\ref{eq16}) ensures the validity of (\ref{sed1}). The verification of \eqref{sed1} for $\langle\oy,v\rangle<0$ is similar.

{\bf Case~3:} $\ox\in\bd(\Q)\setminus\{0\}$. Defining the mapping  $\psi\colon\mathbb{R}^{m+1}\to\mathbb{R}^2$ by
\begin{equation}\label{mappsi}
\psi(x_0,x_r):=\big(\Vert x_r\Vert^2-x_0^2,-x_0\big),\quad(x_0,x_r)\in\mathbb{R}\times\R^{m},
\end{equation}
observe the following representations of the Lorentz cone and its indicator function, respectively:
\begin{equation}\label{lk02}
\Q=\big\{x\in\mathbb{R}^{m+1}\;\big|\;\psi(x)\in\mathbb{R}_-^2\big\}\quad\quad \mbox{and}\quad\quad\delta_\Q=\delta_{\mathbb{R}_-^2}\circ\psi.
\end{equation}
For any $v\in\mathbb{R}^{m+1}$ and $t>0$ we form the vector
\begin{equation}\label{eq14}
w:=\dfrac{\psi(\ox+tv)-\psi(\ox)}{t}
\end{equation}
and use it to write down the relationships
\begin{equation}\label{eq13}
\delta_\Q(\ox+tv)=\delta_{\mathbb{R}_-^2}(\psi(\ox)+tw)\quad\textrm{and}\quad\delta_\Q(\ox)=\delta_{\mathbb{R}_-^2}\big(\psi(\ox)\big).
\end{equation}
It is easy to see that $\nabla\psi(\ox)$ is surjective due to $\ox\in\bd(\Q)\setminus\{0\}$. Employing the first-order chain rule from convex analysis, we get
$N_\Q(\ox)=\nabla\psi(\ox)^*N_{\mathbb{R}_-^2}(\psi(\ox))$. This together with $\oy\in N_\Q(\ox)$ yields the existence of some $\bar\lambda\in N_{\mathbb{R}_-^2}(\psi(\ox))$
for which $\oy=\nabla\psi(\ox)^*\bar\lambda$. This allows us to arrive at
\begin{equation*}
-t\langle\oy,v\rangle=-t\langle\nabla\psi(\ox)^*\bar\lambda,v\rangle=-t\langle\bar\lambda,w\rangle+\langle\bar\lambda,t(w-\nabla\psi(\ox)v)\rangle.
\end{equation*}
Furthermore, it follows from \eqref{eq14} that
\begin{equation*}
t\big(w-\nabla\psi(\ox)v\big)=\psi(\ox+tv)-\psi(\ox)-t\nabla\psi(\ox)v
=\dfrac{1}{2}t^2\big\langle\nabla^2\psi(\ox)v,v\big\rangle+o(t^2),
\end{equation*}
which in turn leads us to the representation
\begin{equation*}
-t\langle\oy,v\rangle=-t\langle\bar\lambda,w\rangle+\dfrac{1}{2}t^2\left\langle\nabla^2\langle\bar\lambda,\psi\rangle(\ox)v,v\right\rangle+o(t^2).
\end{equation*}
Combining the latter with \eqref{eq13} and (\ref{lk01}) readily yields
\begin{eqnarray}\label{eq17}
\Delta_t^2\delta_\Q(\ox|\oy)(v)&=&\dfrac{\delta_{\mathbb{R}_-^2}(\psi(\ox)+t w)-\delta_{\mathbb{R}_-^2}(\psi(\ox))
-t\langle\bar\lambda,w\rangle}{\frac{1}{2}t^2}+\left\langle\nabla^2\langle\bar\lambda,\psi\rangle(\ox)v,v\right\rangle+\dfrac{o(t^2)}{t^2}\nonumber\\
&=&\Delta_t^2\delta_{\mathbb{R}_-^2}\big(\psi(\ox)|\bar\lambda\big)(w)+\left\langle\nabla^2\langle\bar\lambda,\psi\rangle(\ox)v,v\right\rangle+\dfrac{o(t^2)}{t^2}.
\end{eqnarray}
Pick next arbitrary sequences $v_k\to v$ and $t_k\downarrow 0$, and define $w_k:=\dfrac{\psi(\ox+t_kv_k)-\psi(\ox)}{t_k}$ similarly to \eqref{eq14}.
Since $w_k\to\nabla\psi(\ox)v$ as $k\to\infty$, we conclude from \eqref{eq17} that
\begin{eqnarray*}
\liminf_{k\to\infty}\Delta_{t_k}^2\delta_\Q(\ox|\oy)(v_k)&=&\liminf_{k\to\infty}\Big\{\Delta_{t_k}^2\delta_{\mathbb{R}_-^2}\big(\psi(\ox)|\bar\lambda\big)(w_k)
+\left\langle\nabla^2\langle\bar\lambda,\psi\rangle(\ox)v_k,v_k\right\rangle+\dfrac{o(t_k^2)}{t_k^2}\Big\}\\
&\ge&\left\langle\nabla^2\langle\bar\lambda,\psi\rangle(\ox)v,v\right\rangle+\inf_{\tilde w_k\to\nabla\psi(\ox)v}\liminf_{k\to\infty}
\Delta_{t_k}^2\delta_{\mathbb{R}_-^2}\big(\psi(\ox)|\bar\lambda\big)(\tilde w_k)\\
&\ge&\begin{cases}
\left\langle\nabla^2\langle\bar\lambda,\psi\rangle(\ox)v,v\right\rangle&\textrm{if }\;\nabla\psi(\ox)v\in T_{\mathbb{R}_-^2}\big(\psi(\ox)\big)\cap\{\bar\lambda\}^\perp,\\
\infty&\textrm{otherwise},
\end{cases}
\end{eqnarray*}
where the last inequality comes from \cite[Proposition~13.9]{rw} in which the twice epi-differentiability of the indicator function of a convex polyhedron
was established. On the other hand, it follows from the surjectivity of $\nabla\psi(\ox)$ and (\ref{lk02}) that
$$
v\in T_\Q(\ox)\cap\{\oy\}^\perp\Longleftrightarrow\nabla\psi(\ox)v\in T_{\mathbb{R}_-^2}\big(\psi(\ox)\big)\cap\{\bar\lambda\}^\perp,
$$
which in turn leads us to the estimate
\begin{equation}\label{eq21}
\liminf_{k\to\infty}\Delta_{t_k}^2\delta_\Q(\ox|\oy)(v_k)
\ge\begin{cases}
\big\langle\nabla^2\langle\bar\lambda,\psi\rangle(\ox)v,v\big\rangle
&\textrm{if }\;v\in T_\Q(\ox)\cap\{\oy\}^\perp,\\
\infty&\textrm{otherwise}.
\end{cases}
\end{equation}
To finish the proof of (\ref{sed1}), recall that $\bar\lambda\in N_{\mathbb{R}_-^2}(\psi(\ox))$ with $\ox=(\ox_0,\ox_r)\in\bd(\Q)\setminus\{0\}$. Thus we get $\bar\lambda=(\bar\alpha,0)$ for some $\bar\alpha\ge 0$ and so deduce from here and the notation $\widehat{\ox}$ introduced in Section~1 the following equalities:
\begin{equation*}
\oy=\nabla\psi(\ox)^*\bar\lambda=\begin{bmatrix}
-2\ox_0&-1\\
2\ox_r&0\\
\end{bmatrix}\begin{pmatrix}
\bar\alpha\\
0\\
\end{pmatrix}=2 \bar\alpha\widehat{\ox},
\end{equation*}
which yield $\bar\alpha=\dfrac{\Vert\oy\Vert}{2\Vert\widehat{\ox}\Vert}=\dfrac{\Vert\oy\Vert}{2\Vert\ox\Vert}$. Employing now (\ref{mappsi}) brings us to the relationships
\begin{equation*}
\langle\bar\lambda,\psi\rangle(\ox)=\bar\alpha(-\ox_0^2+\Vert\ox_r\Vert^2),\quad\nabla^2\langle\bar\lambda,\psi\rangle(\ox)=2\bar\alpha\begin{bmatrix}
-1&0\\
0&I\\
\end{bmatrix},
\end{equation*}
\begin{equation}\label{eq22}
\big\langle\nabla^2\langle\bar\lambda,\psi\rangle(\ox)v,v\big\rangle=2\bar\alpha(-v_0^2+\Vert v_r\Vert^2)=\dfrac{\Vert\oy\Vert}{\Vert\ox\Vert}(-v_0^2+\Vert v_r\Vert^2).
\end{equation}
Unifying it with \eqref{eq21} verifies the first condition (\ref{sed1}) in the second-order epi-differentiability.

It remains to prove the other condition (\ref{sed2}) in the framework of Case~3. The latter inequality clearly holds when the right-hand side of it equals infinity. Thus we only need
to consider the situation where $v\in{\Bar{\cal K}}$ with the critical cone ${\Bar{\cal K}}$ described by
$$
{\Bar{\cal K}}=T_\Q(\ox)\cap\{\oy\}^\perp=\left\{\begin{array}{ll}
\big\{u\in\mathbb{R}^{m+1}\;\big|\;\langle u,\widehat{\ox}\rangle\le 0\big\}&\mbox{if}\;\oy=0,\\
\big\{u\in\mathbb{R}^{m+1}\;\big|\;\langle u,\widehat{\ox}\rangle=0\big\}&\mbox{if}\;\oy\ne 0.
\end{array}\right.
$$
Construct a sequence $v_k\to v$ satisfying (\ref{sed2}) based on the position of $v$ in ${\Bar{\cal K}}$ as follows:

{\bf Case~3(i):} $v\in\bd({\Bar{\cal K}})\cap\Q$ or $v\in\int({\Bar{\cal K}})$. Having $v=(v_0,v_r)\in \R\times\R^m$, define $v_k:=v$ for any $k$ and claim that
$\ox+tv=(\ox_0+tv_0,\ox_r+t v_r)\in\Q$ when $t>0$ is small enough. This is clear if $v\in\bd({\Bar{\cal K}})\cap\Q$. To justify the claim, it suffices to show that
\begin{equation}\label{lk03}
\ox_0+tv_0 \ge\|\ox_r+tv_r\|
\end{equation}
for all small $t>0$ provided that $v\in\int({\Bar{\cal K}})$. We easily derive that $\langle\widehat{\ox},v\rangle<0$ and $\|\ox_r\|=\ox_0>0$ from the facts that $v\in\int ({\Bar{\cal K}})$ and $\ox=(\ox_0,\ox_r)\in\bd(\Q)\setminus\{0\}$, respectively. This yields
$$
\ox_0+tv_0>0\quad\quad\mbox{and}\quad\quad\la v_r,\ox_r\ra-\ox_0 v_0+t (\|v_r\|^2-v_0^2)<0
$$
for $t$ sufficiently small. The above inequalities tell us that $(\ox_0+tv_0)^2>\|\ox_r+tv_r\|^2$, which thus verifies (\ref{lk03}). Letting $t_k\downarrow 0$, we deduce from $\ox+t_kv\in\Q$ and $v\in\{\oy\}^\perp$ that
\begin{equation}\label{lk04}
\Delta_{t_k}^2\delta_\Q(\ox|\oy)(v_k)=\dfrac{\delta_\Q(\ox+t_k v)-\delta_\Q(\ox)-t_k\langle\oy,v\rangle}{\frac{1}{2}t_k^2}=0
\end{equation}
for $k$ sufficiently large. It is not hard to see furthermore that
$$
\Delta(\ox,\oy)(v)=\dfrac{\Vert\oy\Vert}{\Vert\ox\Vert}(-v_0^2+\Vert v_r\Vert^2)=0.
$$
Combining this with (\ref{lk04}) justifies (\ref{sed2}) under the imposed conditions on $v$.

{\bf Case~3(ii):} $v=(v_0,v_r)\in\bd({\Bar{\cal K}})\setminus Q$. Assume without loss of generality that $\Vert\ox\Vert=\Vert v\Vert=1$. Remembering that $\widehat{\ox}=(-\ox_0,\ox_r)$ according to the notation of Section~1, we conclude from $-\ox_0 v_0+\langle\ox_r,v_r\rangle=\langle\widehat{\ox},v\rangle=0$ and
$\ox=(\ox_0,\ox_r)\in\bd(\Q)\setminus\{0\}$ that
\begin{equation}\label{lk05}
 \|v_r\|^2-v_0^2\ge 0.
\end{equation}
Letting $t_k\downarrow 0$ and employing \eqref{eq17} and \eqref{eq22} yield
\begin{equation}\label{eq24}
\limsup_{k\to\infty}\Delta_{t_k}^2\delta_\Q(\ox|\oy)(v_k)=\limsup_{k\to\infty}\Delta_{t_k}^2\delta_{\mathbb{R}_-^2}\big(\psi(\ox)|\bar \lambda\big)(w_k)
+\dfrac{\Vert\oy\Vert}{\Vert\ox\Vert}(-v_0^2+\Vert v_r\Vert^2).
\end{equation}
Define further the sequence of vectors $v_k$ by
\begin{equation}\label{eq23}
v_k:=\dfrac{x_k-\ox}{t_k}\quad\quad\mbox{with}\quad\quad x_k:=\ox+\alpha_kv-\beta_k\widehat{\ox}\quad\mbox{and}\quad\beta_k=\dfrac{\alpha_k^2(-v_0^2+\Vert v_r\Vert^2)}{4\ox_0(\ox_0+\alpha_kv_0)},
\end{equation}
where $\al_k>0$ is chosen---we will show in the {\em claim} below that such a number $\al_k$ does exist for each $k$---so that $\|x_k-\ox\|=t_k$ and $x_k\in \bd(\Q)$. It follows from construction (\ref{eq23}) of $v_k=(v_{k,0},v_{k,r})\in \R\times\R^m$ that the vectors $w_k$ defined above admit the representations
\begin{equation*}
w_k=\dfrac{\psi(\ox+t_kv_k)-\psi(\ox)}{t_k}=\dfrac{1}{t_k}\Big((0,-\ox_0-t_kv_{k,0})-(0,-\ox_0)\Big)=(0,-v_{k,0}),
\end{equation*}
This tells us that $\langle\bar\lambda,w_k\rangle=\langle(\bar\alpha,0),(0,-v_{k,0})\ra=0$ and implies in turn that
\begin{equation*}
\Delta_{t_k}^2\delta_{\mathbb{R}_-^2}\big(\psi(\ox)|\bar\lambda\big)(w_k)=\dfrac{\delta_{\mathbb{R}_-^2}\big(\psi(\ox+t_kv_k)\big)-\delta_{\mathbb{R}_-^2}\big(\psi(\ox)\big)-t_k\langle \bar\lambda,w_k\rangle}{\frac{1}{2}t_k^2}=0\quad\mbox{for all}\;\;k\in\N.
\end{equation*}
It allows us to arrive at the equality
$$
\limsup_{k\to\infty}\Delta_{t_k}^2\delta_{\mathbb{R}_-^2}\big(\psi(\ox)|\bar\lambda\big)(w_k)=0,
$$
which together with \eqref{eq24} justifies the second twice epi-differentiability requirement (\ref{sed2}).

Let us now verify the aforementioned claim formulated as follows.

{\bf Claim.} {\em For any $v_0\ge 0$ in Case~{\rm 3(ii)} and any $k\in\N$ there is $\al_k>0$ satisfying \eqref{eq23} such that $x_k\in\bd(\Q)$ and $\|x_k-\ox\|=t_k$. If $v_0<0$ in this case, then we can select $\al_k\in(0,-\frac{\ox_0}{v_0})$  as $k\in\N$ so that the above conditions on $x_k$ from \eqref{eq23} are also satisfied}.

We prove this claim by arguing in parallel for both cases of $v_0\ge 0$ and $v_0<0$. Pick $v_0\ge 0$ (resp.\ $v_0<0$) satisfying (\ref{lk05}) and observe that $\beta_k\ge0$ when $\al_k>0$ (resp.\ when $\al_k\in(0,-\frac{\ox_0}{v_0})$) in \eqref{eq23}. Taking into account that $\ox_0^2=\Vert\ox_r\Vert^2$ and that $\ox_0 v_0=\langle\ox_r,v_r\rangle$, we obtain by the direct calculation that the relationship
\begin{equation*}
-\big((1+\beta_k)\ox_0+\alpha_kv_0\big)^2+\left\Vert(1-\beta_k)\ox_r+\alpha_kv_r\right\Vert^2=0
\end{equation*}
is valid in both cases and yields in turn the inequality
$$
\left\Vert(1-\beta_k)\ox_r+\alpha_kv_r\right\Vert=(1+\beta_k)\ox_0+\alpha_k v_0>0.
$$
This confirms that if $v_0\ge 0$ (resp.\ $v_0<0$), then for any $\al_k>0$ (resp.\ $\al_k\in (0,-\frac{\ox_0}{v_0})$) we have
$$
x_k=\big((1+\beta_k)\ox_0+\alpha_kv_0,(1-\beta_k)\ox_r+\alpha_k v_r\big)\in\bd(\Q).
$$
To furnish the verification of the claim, it remains to show that for each $k\in\N$ there exists $\al_k$ from the intervals above such that $\Vert x_k-\ox\Vert=t_k$. To proceed, consider the polynomial
\begin{equation*}
p(\alpha)=\left((-v_0^2+\Vert v_r\Vert^2)^2+16\ox_0^2v_0^2\right)\alpha^4+32\ox_0^3v_0\alpha^3+16(\ox_0^4-t_k^2\ox_0^2v_0^2)\alpha^2-32t_k^2\ox_0^3v_0\alpha-16t_k^2\ox_0^4.
\end{equation*}
Since $p(0)=-16t_k^2\ox_0^4<0$ and the leading coefficient of $p(\al)$ is positive, this polynomial has a positive zero, which we denote by $\al_k$. It follows from
\begin{equation}\label{lk06}
t_k^2=\Vert x_k-\ox\Vert^2=\|\alpha_kv-\beta_k\widehat{\ox}\|^2=\al_k^2+\beta_k^2=\alpha_k^2+\dfrac{\alpha_k^4(-v_0^2+\Vert v_r\Vert^2)^2}{16\ox_0^2(\ox_0+\alpha_kv_0)^2}
\end{equation}
that any root $\al_k>0$ satisfies all our requirements in \eqref{eq23} provided that $v_0\ge 0$. If $v_0<0$, we need to show in addition that there is a root of $p(\al)$ belonging to the interval $(0,-\dfrac{\ox_0}{v_0})$. But it is an immediate consequence of the conditions
\begin{equation*}
p\big (-\dfrac{\ox_0}{v_0}\big)=\dfrac{(-v_0^2+\Vert v_r\Vert^2)^2\ox_0^4}{v_0^4}>0\quad\quad\mbox{and}\quad\quad p(0)=-16t_k^2\ox_0^4<0,
\end{equation*}
which therefore finish the proof of this claim.\vspace*{0.03in}

Let us finally show that $v_k\to v$ as $k\to\infty$. From (\ref{lk06}) we get that $\al_k\to 0$ since $t_k\downarrow 0$ as $k\to\infty$. Remembering that
$\Vert v_k\Vert=1=\Vert v\Vert$, it follows directly from (\ref{eq23}) and (\ref{lk06}) that
$$
\|v_k-v\|^2= 2-2\la v_k,v\ra=2-\frac{2\al_k}{t_k}=2-\frac{2\al_k}{\sqrt{\al_k^2+\beta_k^2}}=2-\frac{2}{\sqrt{1+\frac{\beta_k^2}{\al_k^2}}}\to 2-2=0
$$
as $k\to\infty$, and hence $v_k\to v$. The the proof of the theorem is complete. $\h$\vspace*{0.05in}

In the rest of this section we present some immediate consequences of Theorem~\ref{ted} important in second-order variational analysis of SOCPs. The first one uses the established twice epi-differentiability of $\dd_\Q$ to verify a derivative-coderivative relationship for the normal cone to $\Q$.\vspace*{-0.1in}

\begin{Corollary}{\bf(derivative-coderivative relationship between the normal cone to $\Q$).}\label{dcod2} Let $\ox\in\Q$ and $\oy\in N_\Q(\ox)$. Then we have the inclusion
$$
(DN_\Q)(\ox,\oy)(v)\subset(D^*N_\Q)(\ox,\oy)(v)\quad\quad\mbox{for all}\;\;v\in\R^{m+1}.
$$
\end{Corollary}\vspace*{-0.1in}
{\bf Proof.} It follows from \cite[Theorem~13.57]{rw} that the claimed inclusion holds for any convex set whose indicator function is twice epi-differentiable at the reference point. The latter is the case for the second-order cone $\Q$ due to Theorem~\ref{ted}.$\h$\vspace*{0.05in}

The next corollary provides a precise calculation for the graphical derivative \eqref{gder} of the normal cone to $\Q$ that is significant for the subsequent material of the paper.\vspace*{-0.1in}

\begin{Corollary}{\bf(graphical derivative of the normal cone to $\Q$).}\label{dcod} Let $\ox\in\Q$ and $\oy\in N_\Q(\ox)$. Then for all $v=(v_0,v_r)\in{\Bar{\cal K}}$ the graphical derivative of $N_\Q$ admits the representation
\begin{equation}\label{grder}
(DN_\Q)(\ox,\oy)(v)=\left\{\begin{array}{ll}
N_{\Bar{\cal K}}(v)\quad&\mbox{if }\;\ox\in[\int(\Q)\cup\{0\}],\\
\disp\frac{\|\oy\|}{\|\ox\|}(-v_0,v_r)+N_{\Bar{\cal K}}(v)\quad&\mbox{if }\;\ox\in\bd(\Q)\setminus\{0\},
\end{array}
\right.
\end{equation}
where the critical cone $\Bar{\cal K}$ is defined in Theorem~{\rm\ref{ted}}.
\end{Corollary}\vspace*{-0.1in}
{\bf Proof.} It follows from \cite[Theorem~13.40]{rw} by the twice epi-differentiability of $\dd_\Q$ proved in Theorem~\ref{ted} that for all $v\in\R^{m+1}$ we have
$$
(DN_\Q)(\ox,\oy)(v)=\sub\Big(\frac{1}{2}d^2\delta_\Q(\ox|\oy)\Big)(v).
$$
Combining this with the second epi-derivative formula from Theorem~\ref{ted} verifies the claimed representation of the graphical derivative of the normal cone mapping $x\mapsto N_\Q(x)$. $\h$\vspace*{0.05in}

Now we discuss relationships between the obtained results and a major condition introduced and employed in \cite{mor} for representing the graphical derivative of the normal cone mappings in conic programming under the nondegeneracy condition. Given a closed set $\O\subset\R^n$, assume that the projection operator $\Pi_\O\colon\R^n\tto\R^n$ admits the classical directional derivative $\Pi'_\O(x;h)$ at each $x\in\R^n$ $x$ in any direction $h$. Following \cite[Definition~4.1]{mor}, recall that $\O$ satisfies the {\em projection derivation condition} (PDC) at $x\in\O$ if we have the representation
$$
\Pi'_\O(x+y;h)=\Pi_{{\cal K}(x,y)}(h)\quad\quad\mbox{for all}\;\;y\in N_\O(x)\quad\mbox{and}\;\;h\in\R^n
$$
via the critical cone ${\cal K}(x,y)=T_\O(x)\cap\{y\}^\bot$. It is proved in \cite{mor} that PDC is valid for any convex set $\O$ satisfying the extended polyhedrality condition from \cite[Definition~3.52]{bs} (this includes convex polyhedra) and may also hold in nonpolyhedral settings. Furthermore, PDC holds at the vertex of any convex cone $\O$. On the other hand, we show below that PDC {\em fails} at every nonzero boundary point of the nonpolyhedral Lorentz cone $\Q$ despite its second-order regularity \cite{bs} and other nice properties.

To proceed, we first present a useful characterization of PDC important for its own sake.\vspace*{-0.1in}

\begin{Proposition} {\bf(graphical derivative description of the projection derivation condition).}\label{pdc} Let $\O\subset\R^n$ be a convex set. Then PDC holds at $\ox\in\O$ if and only if
\begin{equation}\label{mn04}
(DN_\O)(\ox,\oy)(v)=N_{{{\cal K}}(\ox,\oy)}(v)\quad\quad\mbox{for all}\;\;\oy\in N_\O(\ox)\quad\mbox{and}\;\;v\in\R^n.
\end{equation}
\end{Proposition}\vspace*{-0.1in}
{\bf Proof.} Assuming that PDC holds at $\ox$, take $\oy\in N_\O(\ox)$ and $v\in\R^n$. To verify the inclusion ``$\subset$" in \eqref{mn04}, pick $w\in(DN_\O)(\ox,\oy)(v)$ and get by definition \eqref{gder} that $(v,w)\in T_{\small{\gph N_\O}}(\ox,\oy)$. It follows from \cite[Proposition~6.17]{rw} that
\begin{equation}\label{mn03}
\Pi_\O(x)=\big(I+N_\O\big)^{-1}(x)\quad\quad\mbox{for any}\;\;x\in\R^n.
\end{equation}
Then elementary tangent cone calculus gives us the representation
\begin{equation}\label{mn02}
\begin{array}{lll}
T_{\small{\gph N_\O}}(\ox,\oy)&=&\big\{(v,w)\big|\;(v+w,v)\in T_{\small{\gph\Pi_\O}}(\ox+\oy,\oy)\big\}\\
&=&\big\{(v,w)\big|\;v=\Pi'_\O(\ox+\oy;v+w)\big\}\;\mbox{ whenever }\;\oy\in N_\O(\ox).
\end{array}
\end{equation}
The above relationships readily imply that
$$
v=\Pi'_\O(\ox+\oy;v+w)=\Pi_{{\cal K}(\ox,\oy)}(v+w)=\big(I+N_{{\cal K}(\ox,\oy)}\big)^{-1}(v+w).
$$
This leads us in turn to $w\in N_{{{\cal K}}(\ox,\oy)}(v)$ and hence justifies the inclusion ``$\subset$" in \eqref{mn04}. The opposite inclusion can be verified similarly.

Conversely, suppose that equality \eqref{mn04} is satisfied. Pick $h\in\R^n$, $\oy\in N_\O(\ox)$, and $v=\Pi'_\O(\ox+\oy;h)$. Employing (\ref{mn02}) tells us that $(v,h-v)\in T_{\small{\gph N_\O}}(\ox,\oy)$, and hence we get $h-v\in N_{{{\cal K}}(\ox,\oy)}(v)$ due to (\ref{mn04}). Combining the latter with (\ref{mn03}) gives us $v=\Pi_{{{\cal K}}(\ox,\oy)}(h)$, which verifies PDC and thus completes the proof of the proposition. $\h$\vspace*{0.05in}

Now we are ready to demonstrate the aforementioned failure of PDC for the second-order cone $\O=\Q$ under consideration on its entire boundary off the origin.\vspace*{-0.1in}

\begin{Corollary}\label{pdpn}{\bf(failure of PDC for the second-order cone at its nonzero boundary points).} Given $\ox\in\Q$, PDC fails whenever $\ox\in\bd(\Q)\setminus\{0\}$.
\end{Corollary}\vspace*{-0.1in}
{\bf Proof.} Suppose on the contrary that PDC holds at some $\ox\in\bd(\Q)\setminus\{0\}$. Employing the graphical derivative formula from Corollary~\ref{dcod} together with the PDC description of Proposition~\ref{pdc} as $\O=\Q$ and ${\cal K}(\ox,\oy)=\Bar{\cal K}$ shows that
$$
N_{\Bar{\cal K}}(v)=(DN_\Q)(\ox,\oy)(v)=\frac{\|\oy\|}{\|\ox\|}(-v_0,v_r)+N_{\Bar{\cal K}}(v)\quad\quad\mbox{for all}\;\;v=(v_0,v_r)\in\R\times\R^m.
$$
It yields $(-v_0,v_r)\in N_{\Bar{\cal K}}(v)={\Bar{\cal K}}^*\cap\{v\}^\bot$ whenever $\|\oy\|\ne 0$ and $v\in\R^{m+1}$, which is clearly wrong. The obtained contradiction justifies the claimed statement. $\h$\vspace*{-0.1in}

\section{Remarkable Properties of Second-Order Cone Constraints}\label{sec04}\sce\vspace*{-0.05in}

In this section we derive new properties of the second-order cone $\Q$, which are important in what follows while being also of their own interest. The derivation of some of the results below employs those obtained in the previous section.

Our first result here provides a complete description of the set of Lagrange multipliers associated with stationary points of the constraint system $\Gamma$ in \eqref{CS}. Given a stationary pair $(x,x^*)\in\gph N_\Gamma$, define the set of {\em Lagrange multipliers} associated with $(x,x^*)$ by
\begin{equation}\label{lagn}
\Lambda(x,x^*):=\left\{\lambda\in N_\Q\big(\Phi(x)\big)\,\big|\,\nabla\Phi(x)^*\lambda=x^*\right\}
\end{equation}
and the {\em critical cone} to $\Gamma$ at $(x,x^*)$ by
\begin{equation}\label{critco}
K(x,x^*):=T_\Gamma(x)\cap\{x^*\}^\perp.
\end{equation}
If $\Phi(\ox)=0$ for some $\ox$ with $(\bar x,\bar x^*)\in\gph N_\Gamma$, then the Lagrange multiplier set reduces to
\begin{equation}\label{lag0}
\Lambda(\ox,\ox^*)=\left\{\lambda\in-\Q\,\big|\,\nabla\Phi(\ox)^*\lambda=\ox^*\right\}.
\end{equation}

Following \cite[Definition~2.105]{bs}, we say that the {\em Slater condition} holds for $\Lambda(\bar x,\bar x^*)$ if there is a multiplier $\lm\in\int(-\Q)$ such that $\nabla\Phi(\bar x)^*\lm=\ox^*$. The next result provides a precise description of the Lagrange multiplier set \eqref{lag0} that plays a significant role in our method of conducting the second-order analysis of $\Gamma$. A part of this analysis is inspired by the unpublished work of Shapiro and Nemirovski \cite{sn} about the ``no duality gap" property in linear conic programs generated by convex cones; see, in particular, the proof of \cite[Proposition~3]{sn} and the discussion after it.\vspace*{-0.1in}

\begin{Proposition}{\bf(description of Lagrange multipliers for the second-order cone).}\label{lag10} Let $(\bar x,\bar x^*)\in\gph N_\Gamma$ with $\Phi(\ox)=0$, and let $\Lambda(\ox,\ox^*)\ne\emp$ for the set of Lagrange multipliers \eqref{lag0}. Then one of the following alternatives holds for $\Lambda(\ox,\ox^*)$:

{\bf(LMS1)} The set $\Lambda(\bar x,\bar x^*)$ satisfies the Slater condition, i.e., there exists a Lagrange multiplier $\lm\in\int(-\Q)$ with $\nabla\Phi(\bar x)^*\lm=\ox^*$. In this case we get that for any $\olm\in\Lambda(\ox,\ox^*)$ there are numbers $\ell,\ve>0$ ensuring the error bound estimate
\begin{equation}\label{error}
\dist\big(\lm;\Lambda(\ox,\ox^*)\big)\le\ell\big(\dist(\lm;-\Q)+\|\nabla\Phi(\ox)^*\lm-\ox^*\|\big)\quad\mbox{whenever}\;\;\lm\in\B_\ve(\olm).
\end{equation}

{\bf(LMS2)} $\Lambda(\ox,\ox^*)=\{\olm\}$ for some multiplier $\olm\in\bd{(-\Q)}\setminus\{0\}$.

{\bf(LMS3)} $\Lambda(\ox,\ox^*)=\big\{t\olm\,\big|\,t\ge 0\big\}$ for some $\olm\in\bd{(-\Q)}$. In this case we have $\ox^*=0$.
\end{Proposition}\vspace*{-0.1in}
{\bf Proof.} The validity of the error bound estimate \eqref{error} in the Slater case (LMS1) follows from \cite[Corollary~5]{bbl}. Suppose that the Slater condition fails and pick any $\olm\in\Lambda(\ox,\ox^*)$. This obviously ensures the fulfillments of either (LMS2) or (LMS3) provided that $\Lambda(\ox,\ox^*)$ is a singleton. If the latter doesn't hold, we claim that $\Lambda(\ox,\ox^*)\subset\R_+\olm$. Assuming the contrary yields $\olm\ne 0$ and allows us to find $0\ne\lm\in\Lambda(\ox,\ox^*)$ such that $\lm\not\in\R_+\olm$. Since the Slater condition fails, we have $\olm,\,\lm\in\bd(-\Q)\setminus\{0\}$. Define now $\lm_\al:=\al\olm+(1-\al)\lm$ with $\al\in(0,1)$ and observe that $\lm_\al\in\int(-\Q)$; otherwise we get $\lm\in\R_+\olm$. This observation amounts to saying that the Slater condition holds for $\Lambda(\ox,\ox^*)$, which is a contradiction. Thus we arrive at the inclusion $\Lambda(\ox,\ox^*)\subset\R_+\olm$ telling us that either (LMS2) or (LMS3) is satisfied. Since $0\in\Lambda(\ox,\ox^*)$ in case (LMS3), we get $\ox^*=0$ in this case and hence complete the proof of the proposition. $\h$\vspace*{0.05in}

To proceed with our further analysis, we introduce an appropriate (very weak) constraint qualification for the second-order cone constraint system \eqref{CS}. This condition has been recently employed in the polyhedral framework of NLPs to conduct a second-order analysis of the classical equality and inequality constraint systems with ${\cal C}^2$-smooth data; see \cite{ch,gm15,go16}. It has also been studied in \cite{gm16} in nonpolyhedral settings via first-order and second-order constructions of variational analysis. However, to the best of our knowledge, it has never been implemented before for the second-order variational analysis of nonpolyhedral systems as we do in this paper.\vspace*{-0.1in}

\begin{Definition}{\bf(metric subregularity constraint qualification).}\label{defmscq} We say that system \eqref{CS}  satisfies the {\sc metric subregularity constraint qualification} $($MSCQ$)$ at $\bar x\in\Gamma$ with modulus $\kappa>0$ if the mapping $x\mapsto\Phi(x)-\Q$ is metrically subregular at $(\bar x,0)$ with modulus $\kappa$.
\end{Definition}\vspace*{-0.1in}

Using \eqref{metreq} with the fixed vector $y=\oy=0$, observe that the introduced MSCQ with modulus $\kappa$ for the constraint system \eqref{CS} can be equivalently described as the existence of a neighborhood $U$ of $\ox$ such that the distance estimate
\begin{equation}\label{mscq}
\dist(x;\Gamma)\le\kappa\,\dist\big(\Phi(x);\Q\big)\quad\quad\mbox{for all}\;\;x\in U
\end{equation}
holds. It is worth mentioning that the defined MSCQ property of \eqref{CS} is {\em robust} in the sense that its validity at $\ox\in\Gamma$ yields this property at any $x\in\Gamma$ near $\ox$. Furthermore, it is clear Example~\ref{Ex1} below) that the MSCQ from Definition~\ref{defmscq} is strictly weaker that the qualification condition corresponding to the {\em metric regularity} of the mapping $x\mapsto\Phi(x)-\Q$ around $(\ox,0)$ therein. The latter is well known to be equivalent to the {\em Robinson constraint qualification} (RCQ), which is the basic qualification condition in conic programming:
\begin{equation}\label{gf07}
N_\Q\big(\Phi(\ox)\big)\cap\ker\nabla\Phi(\ox)^*=\{0\}.
\end{equation}
An important role of MSCQ and its calmness equivalent for inverse mappings has been recognized in generalized differential calculus of variational analysis. In particular, it follows from \cite[Theorem~4.1]{hjo} and the convexity of $\Q$ that there is a neighborhood $U$ of $\ox$ such that
\begin{equation}\label{first}
N_\Gamma(x)=\widehat{N}_\Gamma(x)=\nabla\Phi(x)^*N_\Q\big(\Phi(x)\big)\quad\quad\mbox{for all}\;\;x\in\Gamma\cap U,
\end{equation}
where $\Hat N_\O(\ox)$ stands for the {\em regular/Fr\'echet normal cone} to $\O$ at $\ox\in\O$ defined by
\begin{eqnarray}\label{rnc}
\Hat N_\O(\ox):=\Big\{v\in\R^n\Big|\;\limsup_{x\st{\O}{\to}\ox}\frac{\la v,x-\ox\ra}{\|x-\ox\|}\le 0\Big\},
\end{eqnarray}
which is {\em dual} to the tangent cone \eqref{2.5}, i.e, $\Hat N_\O(\ox)=T^*_\O(\ox)$. The first equality in \eqref{first} postulates the {\em normal regularity} of $\Gamma$ at any point $x\in\Gamma$ near $\ox$. Note also that the validity of MSCQ for $\Gamma$ at $\ox\in\Gamma$ ensures by \cite[Proposition~1]{ho} the tangent cone calculus rule
\begin{equation}\label{first2}
T_\Gamma(x)=\big\{v\in\R^{n}\big|\;\nabla\Phi(x)v\in T_\Q\big(\Phi(x)\big)\big\}\quad\quad\mbox{for all}\;\;x\in\Gamma\cap U.
\end{equation}

To proceed further, recall that the second-order cone $\Q$ is {\em reducible} at its {\em nonzero boundary points} to a convex polyhedron in the sense of \cite[Definition~3.135]{bs}; this was first shown in \cite[Lemma~15]{BoR05}. In what follows we use a different reduction of $\Q$ via the mapping $\psi$ from \eqref{mappsi} that allows us to simplify the subsequent calculations. Indeed, the alternative representation \eqref{lk02} of the second-order cone $\Q$ via the mapping $\psi$ from \eqref{mappsi} in the proof of Case~3 of Theorem~\ref{ted} is instrumental to furnish the reduction of $\Q$ to $\R^2_-$ at nonzero boundary points $x$. Observe that the Jacobian matrix $\nabla\psi(x)$ has full rank and get the representation
\begin{equation}\label{CS1}
\Gamma=\left\{x\in\mathbb{R}^n\,\big|\,(\psi\circ\Phi)(x)\in\mathbb{R}_-^2\right\}\quad\quad\mbox{whenever}\;\;\Phi(x)\in\bd(\Q)\setminus\{0\}.
\end{equation}

By showing below that the metric subregularity of the mapping $x\mapsto\Phi(x)-\Q$ at nonzero boundary points yields the one for $x\mapsto(\psi\circ\Phi)(x)-\R^2_-$, we open the door to the usage in this case the results for convex polyhedra established \cite{go16}. It is convenient to implement the decomposition of the vectors $\Phi(x)\in\R^{m+1}$ relevant to that in the second-order cone \eqref{Qm}:
\begin{equation}\label{phidec}
\Phi(x)=\big(\Phi_0(x),\Phi_r(x)\big)\in\R\times\R^m\quad\quad\mbox{as}\;\;x\in\R^n.
\end{equation}\vspace*{-0.3in}

\begin{Lemma}{\bf(propagation of metric subregularity for nonzero boundary points of the second-order cone).}\label{mscq2} Let $\bar x\in\Gamma$ be such that $\Phi(\ox)\in\bd(\Q)\setminus\{0\}$. Then the metric subregularity of the mapping $x\mapsto\Phi(x)-\Q$ at $(\bar x,0)$ ensures the one for
$x\mapsto(\psi\circ\Phi)(x)-\R^2_-$ at $(\bar x,0)$ with the mapping $\psi\colon\R^{m+1}\to\R^2$ taken from \eqref{mappsi}.\vspace*{-0.1in}
\end{Lemma}
{\bf Proof.} To verify the lemma, we need to establish the existence of a positive number $\kappa$ and a neighborhood $V$ of $\ox$ such that the estimate
\begin{equation}\label{mscq3}
\dist(x;\Gamma)\le\kappa\,\dist\big((\psi\circ\Phi)(x);\mathbb{R}_-^2\big)\quad\quad\mbox{for all}\;\;x\in V
\end{equation}
holds. Let us first show that there are a constant $c>0$ and a neighborhood $U$ of $\bar x$ for which
\begin{equation}\label{eq1}
\dist\big(\Phi(x);\Q\big)\le c\;\dist\big((\psi\circ\Phi)(x);\mathbb{R}_-^2\big)\quad\quad\mbox{for all}\;\;x\in{U}.
\end{equation}
Indeed, employing \eqref{phidec} together with the direct calculations tells us that
\begin{equation}\label{bg01}
\dist\big(\Phi(x);\Q\big)=\begin{cases}
0&\textrm{if }\;\Phi(x)\in\Q,\\
\Vert\Phi(x)\Vert&\textrm{if }\;\Phi(x)\in-\Q,\\
\dfrac{\sqrt{2}}{2}\big(\Vert\Phi_r(x)\Vert-\Phi_0(x)\big)&\textrm{if }\;\Phi(x)\notin\Q\cup(-\Q);
\end{cases}
\end{equation}
\begin{equation*}
\dist\big((\psi\circ\Phi)(x);\mathbb{R}_-^2\big)=\begin{cases}
0&\textrm{if }\;\Phi(x)\in\Q,\\
-\Phi_0(x)&\textrm{if }\Phi(x)\in-\Q,\\
\Vert \Phi_r(x)\Vert^2-\Phi_0^2(x)&\textrm{if }\;\Phi(x)\notin\Q\cup(-\Q)\;\;\mbox{and}\;\;\Phi_0(x)\ge 0,\\
\sqrt{(\Vert\Phi_r(x)\Vert^2-\Phi_0^2(x))^2+\Phi_0^2(x)}&\textrm{if }\;\Phi(x)\notin\Q\cup(-\Q)\;\;\mbox{and}\;\;\Phi_0(x)<0.
\end{cases}
\end{equation*}
It follows from $\bar x\in\Gamma$ and $\Phi_0(\bar x)=\|\Phi_r(\ox)\|\ne 0$ that there exists a neighborhood $U$ of $\bar x$ such that the inequality $\Phi_0(x)>\dfrac{1}{2}\Phi_0(\bar x)$ holds whenever $x\in U$. Pick $x\in U$ and observe that the two cases may occur: either (a) $\Phi(x)\in\Q$ for which we have  $\dist(\Phi(x);\Q)=\dist((\psi\circ\Phi)(x);\mathbb{R}_-^2)=0$, and hence estimate \eqref{eq1} is clearly satisfied, or (b) $\Phi(x)\notin\Q$, which means by \eqref{phidec} that $\Vert\Phi_r(x)\Vert>\Phi_0(x)$. Therefore we arrive at
\begin{eqnarray*}
\dist\big((\psi\circ\Phi)(x);\mathbb{R}_-^2\big)
&=&\big(\Vert\Phi_r(x)\Vert-\Phi_0(x)\big)\big(\Vert\Phi_r(x)\Vert+\Phi_0(x)\big)\\
&\ge& 2\sqrt{2}\;\Phi_0(x)\dist\big(\Phi(x);\Q\big)\\
&\ge&\sqrt{2}\;\Phi_0(\bar x)\dist\big(\Phi(x);\Q\big),
\end{eqnarray*}
which justifies estimate \eqref{eq1} with $c:=\big(\sqrt{2}\Phi_0(\bar x)\big)^{-1}$. Combining this and estimate \eqref{mscq} leads us to \eqref{mscq3} and thus completes the proof of the proposition. $\h$\vspace*{0.05in}

The next result is of its own interest while being important for calculating the graphical derivative of the normal cone mapping given in the next section.\vspace*{-0.1in}

\begin{Theorem}{\bf(normal cone to the critical cone of ice-cream constraint systems).}\label{norco} Let $(\bar x,\bar x^*)\in\gph N_\Gamma$. Assuming the validity of MSCQ at $\ox\in\Gamma$ and picking any $\lambda\in\Lambda(\bar x,\bar x^*)$ and $v\in K(\bar x,\bar x^*)$, we represent the normal cone to the critical cone $K(\bar x,\bar x^*)$ by
\begin{equation}\label{norcri}
N_{K(\bar x,\bar x^*)}(v)=\widehat N_{K(\bar x,\bar x^*)}(v)=\nabla\Phi(\bar x)^*\big[T_{N_\Q(\Phi(\bar x))}(\lambda)\cap\{\nabla\Phi(\bar x)v\}^\bot\big].
\end{equation}
\end{Theorem}\vspace*{-0.1in}
{\bf Proof.} It follows from \cite[Corollary~16.4.2]{r70}, \eqref{first}, and the normal-tangent duality that
\begin{equation}\label{eq8}
\big (K(\bar x,\bar x^*)\big)^*=\big (T_\Gamma(\bar x)\cap\{\bar x^*\}^\perp\big)^*={\rm cl}\big(N_\Gamma(\bar x)+\mathbb{R}\bar x^*\big).
\end{equation}
We proceed with verifying the following statement:

{\bf Claim.} {\em If $\Phi(\ox)\in\Q\setminus\{0\}$, then
\begin{equation}\label{ki01}
{\rm cl}\left(N_\Gamma(\bar x)+\mathbb{R}\bar x^*\right)=N_\Gamma(\bar x)+\mathbb{R}\bar x^*.
\end{equation}
Furthermore, \eqref{ki01} is also valid if $\Phi(\ox)=0$ and either {\rm(LMS1)} or {\rm(LMS3)} above holds.}\vspace*{0.03in}

To justify the claim, we split the arguments into the three cases depending on the position of the vector $\Phi(\ox)$ in  the second-order cone $\Q$:

{\bf{Case~1:}} $\Phi(\ox)\in\int(\Q)$. This gives us $\ox^*=0$, which immediately yields \eqref{ki01}.

{\bf{Case~2:}} $\Phi(\ox)\in\bd(\Q)\setminus\{0\}$. Then the normal cone to $\Gamma$ at $\ox$ is a convex polyhedron. Using this together with \cite[Corollary~19.3.2]{r70} ensures the validity of \eqref{ki01}.

{\bf{Case~3:}} $\Phi(\ox)=0$ and either (LMS1) or (LMS3) holds. If the Slater condition in (LMS1) is satisfied, we have $\lm\in\int(-\Q)$ such that $\nabla\Phi(\bar x)^*\lm=\ox^*$, which shows together with \eqref{first} that
$$
N_\Gamma(\bar x)+\mathbb{R}\bar x^*=\nabla\Phi(\bar x)^*N_\Q\big(\Phi(\ox)\big)+\R\nabla\Phi(\bar x)^*\lm =\nabla\Phi(\bar x)^*\big(-\Q +\R\lm\big).
$$
Pick $\eta\in\R^{m+1}$ and find $t>0$ sufficiently small so that $\lm+t\eta\in-\Q$. This leads us to
$$
t\eta=\lm+t\eta-\lm\in-\Q+\R\lm,
$$
and therefore we get $\eta\in -\Q+\R\lm$. It tells us that $-\Q+\R\lm=\R^{m+1}$, which results in
$$
N_\Gamma(\bar x)+\mathbb{R}\bar x^*=\nabla\Phi(\bar x)^*\big(-\Q +\R\lm\big)=\nabla\Phi(\bar x)^*\R^{m+1},
$$
and hence verifies \eqref{ki01} in this setting. To finish the proof of the claim, it remains to recall that under (LMS3) we have $\ox^*=0$, and thus \eqref{ki01} is satisfied.\vspace*{0.03in}

We proceed with the proof of the theorem, we check first that \eqref{norcri} holds for all the cases in the above claim. Picking any $\lambda\in\Lambda(\bar x,\bar x^*)$ and $v\in K(\bar x,\bar x^*)$, deduce from \eqref{ki01} that
\begin{equation}\label{df01}
N_{K(\bar x,\bar x^*)}(v)=\widehat N_{K(\bar x,\bar x^*)}(v)=\big(K(\bar x,\bar x^*)\big)^*\cap\{v\}^\perp=\left (N_\Gamma(\bar x)+\mathbb{R}\bar x^*\right)\cap\{v\}^\perp.
\end{equation}
For each $v^*\in N_{K(\bar x,\bar x^*)}(v)$ we find by \eqref{first} and \eqref{df01} some $\tilde\mu\in N_Q(\bar y)$ and $\alpha\in \mathbb{R}$ with
\begin{equation*}
v^*=\nabla\Phi(\bar x)^*\tilde\mu+\alpha \bar x^*=\nabla\Phi(\bar x)^*(\tilde\mu+\alpha\lambda).
\end{equation*}
Letting $\mu:=\tilde\mu+\alpha\lambda$, we get $\lambda+\varepsilon\mu=(1+\varepsilon\alpha)\lambda+\varepsilon\tilde\mu\in N_Q(\bar y)$ for any small $\varepsilon\geq 0$, which leads us to the inclusion $\mu\in T_{N_Q(\bar y)}(\lambda)$. Taking it into account and using \eqref{df01} give us $\langle\mu,\nabla\Phi(\bar x),v\rangle=\langle v^*,v\rangle=0$, and thus show that $v^*$ belongs to the set on the right-hand side of \eqref{norcri}.

To verify the opposite inclusion in \eqref{norcri}, pick $\mu\in T_{N_Q(\bar y)}(\lambda)$ with $\langle\mu,\nabla\Phi(\bar x)v\rangle=0$ and find sequences $t_k\downarrow 0$ and $\mu_k\to\mu$ with $\lambda+t_k\mu_k\in N_Q(\bar y)$ for all $k\in\N$. It follows from \eqref{first} that
$$
\nabla\Phi(\bar x)^*(\lambda+t_k\mu_k)\in N_\Gamma(\bar x)=\big(T_\Gamma(\bar x)\big)^*.
$$
Using this, for any $w\in K(\bar x,\bar x^*)$ we get
\begin{equation*}
t_k\langle\mu_k,\nabla\Phi(\bar x)w\rangle=\langle\bar x^*,w\rangle+t_k\langle\mu_k,\nabla\Phi(\bar x)w\rangle=\langle\lambda+t_k\mu_k,\nabla\Phi(\bar x)w\rangle\le 0.
\end{equation*}
The passage to the limit as $k\to\infty$ gives us the relationships
$$
\langle\nabla\Phi(\bar x)^*\mu,w\rangle=\langle\mu,\nabla\Phi(\bar x)w\rangle\le 0,
$$
which imply that $\nabla\Phi(\bar x)^*\mu\in\big(K(\bar x,\bar x^*)\big)^*$. Combining it with \eqref{df01} and $\langle\mu,\nabla\Phi(\bar x)v\rangle=0$ leads us to $\nabla\Phi(\bar x)^*\mu\in \widehat N_{K(\bar x,\bar x^*)}(v)$, and thus justifies the inclusion ``$\supset"$ in \eqref{norcri} and the equality therein under the assumptions of the above claim.

Continuing the proof of the theorem, it remains to justify \eqref{norcri} in the setting where $\Phi(\ox)=0$ and (LMS2) holds. Since $\Lambda(\bar x,\bar x^*)=\{\olm\}$ with $\olm=(\olm_0,\olm_r)\in\bd(-\Q)\setminus\{0\}$ in this case, and since MSCQ is satisfied at $\ox$, we have by using \eqref{first2} that
\begin{eqnarray*}
K(\bar x,\bar x^*)&=&T_\Gamma(\bar x)\cap\{\bar x^*\}^\perp=\big\{v\in\mathbb{R}^n\;\big|\;\nabla\Phi(\bar x)v\in\Q\;\;\mbox{and}\;\;\langle v,\nabla\Phi(\bar x)^*\bar\lambda\rangle=0\big\}\\
&=& \big\{v\in\mathbb{R}^n\;\big|\;\nabla\Phi(\bar x)v\in \Q \;\;\mbox{and}\;\;\langle\nabla\Phi(\bar x)v,\bar\lambda\rangle=0\big\}\\
&=&\big\{v\in\mathbb{R}^n\;\big|\;\nabla\Phi(\bar x)v\in\Q\cap\{\bar\lambda\}^\perp\big\}=\big\{v\in\mathbb{R}^n\;\big|\;\nabla\Phi(\bar x)v\in\mathbb{R}_+\hat{\bar\lambda}\big\},
\end{eqnarray*}
where $\hat{\olm}=(-\olm_0,\olm_r)$. Pick now $v\in K(\bar x, \bar x^*)$ and observe that
$$
\begin{array}{ll}
{N}_{K(\bar x,\bar x^*)}(v)&=\nabla\Phi(\ox)^*N_{\R_+\hat{\olm}}\big(\nabla\Phi(\ox)v\big)=\nabla\Phi(\ox)^*\big[\big(\R_+\hat{\olm}\big)^*\cap\{\nabla\Phi(\ox)v\}^\perp\big]\\
&=\nabla\Phi(\ox)^*\big[T_{-\Q}(\olm)\cap\{\nabla\Phi(\ox)v\}^\perp\big]=\nabla\Phi(\ox)^*\big[T_{N_\Q(\oz)}(\olm)\cap\{\nabla\Phi(\ox)v\}^\perp\big],
\end{array}
$$
where the first equality (chain rule) holds by Robinson's seminal result from \cite{rob1} since $\mathbb{R}_+\hat{\bar\lambda}$ is a convex polyhedron and the constraint mapping $\nabla\Phi(\bar x)v$ is linear. This justifies \eqref{norcri} in the case under consideration and thus completes the proof of the theorem.
$\h$\vspace*{0.05in}

A similar result to Theorem~\ref{norco} was established in \cite[Lemma~1]{go16} for polyhedral constraint systems with equality and inequality constraints coming from problems of nonlinear programming. The nonpolyhedral nature of the second-order cone $\Q$ creates significant difficulties in comparison with the polyhedral NLP structure that are successfully overcome in the proof above.\vspace*{0.03in}

Now we present the main result of this section giving a characterization of the simultaneous fulfillment of the {\em uniqueness} of Lagrange multipliers associated with stationary points of \eqref{CS} and a certain {\em error bound} estimate, which is automatic for polyhedral systems. Both properties are algorithmically important; see, e.g., the book \cite{is14} that strongly employs the uniqueness of Lagrange multipliers in polyhedral NLP systems and its characterization via the {\em strict Mangasarian-Fromovitz constraint qualification condition} (SMFCQ) for Newton-type methods.

While dealing with the set $\Gamma$ in the next theorem, the only point $\ox$ that needs to be taken care of is the one for which $\Phi(\ox)=0$. This comes from the observation made right before Lemma~\ref{mscq2} on the reducibility of $\Q$ at nonzero boundary points to the convex polyhedron $\R^2_-$.\vspace*{-0.1in}

\begin{Theorem}{\bf(characterization of uniqueness of Lagrange multipliers with error bound estimate for second-order cone constraints).}\label{unique} Let $(\bar x,\bar x^*)\in\gph N_\Gamma$, and let $\olm\in\Lambda(\ox,\ox^*)$ with $\Phi(\ox)=0$. Then the following statements are equivalent:

{\bf(i)} $\olm$ is a unique multiplier, and for some $\ell>0$ the error bound estimate holds:
\begin{equation}\label{gf01}
\dist(\lm;\Lambda(\ox,\ox^*))\le\ell \;\|\nabla\Phi(\ox)^*\lm-\ox^*\|\quad\quad\mbox{for all}\,\,\;\;\lm\in-\Q.
\end{equation}

{\bf(ii)} The dual qualification condition is satisfied:
\begin{equation}\label{gf02}
(DN_\Q)\big(\Phi(\ox),\olm\big)(0)\cap\ker\nabla\Phi(\ox)^*=\{0\}.
\end{equation}
If in this case $\olm\in\bd(-\Q)\setminus\{0\}$, then \eqref{gf02} implies that the matrix $\nabla\Phi(\ox)$ has full rank.

{\bf(iii)} The strict Robinson constraint qualification holds:
\begin{equation}\label{sqc}
\nabla\Phi(\ox)\R^n-T_\Q\big(\Phi(\ox)\big)\cap\{\olm\}^\bot=\Rm.
\end{equation}
\end{Theorem}\vspace*{-0.1in}
{\bf Proof}. Assume that (ii) is satisfied and pick any $\lm\in\Lambda(\ox,\ox^*)$. We first show that $\lm=\olm$, which verifies the uniqueness of Lagrange multipliers. It readily follows from \eqref{lag0} that
\begin{equation}\label{gf03}
\lm-\olm\in\ker\nabla\Phi(\ox)^*\quad\mbox{and}\quad\lm-\olm\in-\Q+\R\olm.
\end{equation}
Then Corollary~\ref{dcod} tells us that $(DN_\Q)(\oz,\olm)(0)=N_{\Bar{\cal K}}(0)={\Bar{\cal K}}^*$ with ${\Bar{\cal K}}=T_\Q(\Phi(\ox))\cap\{\olm\}^\bot=\Q\cap\{\olm\}^\bot $. Therefore we arrive at the relationships
\begin{equation}\label{gf04}
\lm-\olm\in-\Q+\R\olm\subset\big(\Q\cap\{\olm\}^\bot\big)^*=(DN_\Q)(\Phi(\ox),\olm)(0).
\end{equation}
Using them together with \eqref{gf02} and the first inclusion in \eqref{gf03}, we get $\lm=\olm$.

To verify now the error bound \eqref{gf01} in (i), we use $\Lambda(\ox,\ox^*)=\{\olm\}$ and arguing by contradiction suppose that for any $k\in\N$ there is $\lm_k\in-\Q$ satisfying the conditions
$$
\|\lm_k-\olm\|>k\|\nabla\Phi(\ox)^*\lm_k-\ox^*\|=k\|\nabla\Phi(\ox)^*(\lm_k-\olm)\|.
$$
Assume without loss of generality that $\frac{\lm_k-\olm}{\|\lm_k-\olm\|}\to\eta$ as $k\to \infty$ with $\|\eta\|=1$. Thus  passing to the limit in the above inequality brings us to
\begin{equation}\label{gf05}
\nabla \Phi(\ox)^*\eta=0.
\end{equation}
On the other hand, we have the inclusions
$$
\frac{\lm_k-\olm}{\|\lm_k-\olm\|}\in-\Q+\R\olm\subset\big(\Q\cap\{\olm\}^\bot\big)^*,
$$
which together with \eqref{gf04} ensure the relationships
$$
\eta\in\big(\Q\cap\{\olm\}^\bot\big)^*=(DN_\Q)\big(\Phi(\ox),\olm\big)(0).
$$
Combining the latter with \eqref{gf05} and taking into account (ii) lead us to $\eta=0$, which contradicts the fact that $\|\eta\|=1$ and thus justifies the error bound estimate \eqref{gf01} in (i).

To verify next the converse implication $\rm{(i)}\Longrightarrow\rm{(ii)}$, take $\eta\in(DN_\Q)(\Phi(\ox),\olm)(0)\cap\ker\nabla\Phi(\ox)^*$ and get by the definition of the graphical derivative that $(0,\eta)\in T_{\small{\gph N_\Q}}(\Phi(\ox),\olm)$. This allows us to find sequences $t_k\downarrow 0$ and $(v_k,\eta_k)\to(0,\eta)$ as $k\to\infty$ such that $(\Phi(\ox),\olm)+t_k(v_k,\eta_k)\in\gph N_\Q$ and therefore $\olm+t_k\eta_k\in N_\Q(\Phi(\ox)+t_kv_k)\subset-\Q$. Employing estimate \eqref{gf01} brings us to
$$
\|\olm+t_k\eta_k-\olm\|=\dist\big(\olm+t_k\eta_k;\Lambda(\ox,\ox^*)\big)\le\ell\|\nabla\Phi(\ox)^*(\olm+t_k\eta_k)-\ox^*\|,
$$
which implies in turn that $\|\eta_k\|\le\ell\|\nabla\Phi(\ox)^*\eta_k\|$. Passing to the limit as $k\to\infty$ tells us that $\|\eta\|\le\ell\|\nabla\Phi(\ox)^*\eta\|$. By $\eta\in\ker\nabla\Phi(\ox)^*$ we get $\eta=0$ and thus arrive at \eqref{gf02}.

To finish the proof of (ii), suppose that $\olm=(\olm_0,\olm_r)\in\bd(-\Q)\setminus\{0\}$ and conclude from the graphical derivative formula in Corollary~\ref{dcod} that
$$
(DN_{\Q})\big(\Phi(\ox),\olm\big)(0)=\big(\Q\cap\{\olm\}^\perp\big)^*=\big(\R_+\hat{\olm}\big)^*=\big\{(w_0,w_m)\in\R\times\R^{m}\big|\,\la w_r,\olm_r\ra-w_0\olm_0\le 0\big\}.
$$
It gives us by \eqref{gf02} that ${\rm ker}\nabla\Phi(\ox)^*=\{0\}$, and thus the matrix $\nabla\Phi(\ox)$ is of full rank.

To complete the proof of the theorem, it remains to show that the qualification conditions \eqref{gf02} and \eqref{sqc} are equivalent for the case of \eqref{CS}. Indeed, it follows from \eqref{sqc} that
$$
\big(T_\Q(\Phi(\ox))\cap\{\olm\}^\bot\big)^*\cap\ker\nabla\Phi(\ox)^*=\big(T_\Q(\Phi(\ox))\cap\{\olm\}^\bot-\nabla\Phi(\ox)\R^n\big)^*=\{0\},
$$
and hence the dual qualification condition \eqref{gf02} holds by Corollary~\ref{dcod}. To verify the converse implication, we deduce from \eqref{gf02} that
$$
{\rm cl}\big(\nabla\Phi(\ox)\R^n-T_\Q(\Phi(\ox))\cap\{\olm\}^\bot\big)=\Rm.
$$
Since $\nabla\Phi(\ox)\R^n-T_\Q(\Phi(\ox))\cap\{\olm\}^\bot$ is convex, it has nonempty relative interior. Hence it follows from \cite[Proposition~2.40]{rw} that the relationships
\begin{eqnarray*}
\Rm=\ri(\Rm)&=&\ri\big[{\rm cl}\big(\nabla\Phi(\ox)\R^n-T_\Q\big(\Phi(\ox)\big)\cap\{\olm\}^\bot\big)\big]\\
&=&\ri\big(\nabla\Phi(\ox)\R^n-T_\Q\big(\Phi(\ox)\big)\cap\{\olm\}^\bot\big)\\
&\subset&\big(\nabla\Phi(\ox)\R^n-T_\Q\big(\Phi(\ox)\big)\cap\{\olm\}^\bot\big)
\end{eqnarray*}
are satisfied. This justifies \eqref{sqc} and thus ends the proof of the theorem. $\h$\vspace*{-0.1in}

\begin{Remark}{\bf(discussions on constraint qualifications for second-order cone systems).}\label{gf06}{\rm

{\bf (i)} Condition \eqref{sqc} was introduced in \cite{bs} as ``strict constraint qualification" in conic programming and then was called ``strict Robinson constraint qualification" (SRCQ) in \cite{dsz}. In the case of NLPs this condition reduces to the strict Mangasarian-Fromovitz constraint qualification (SMFCQ) discussed before the formulation of Theorem~\ref{unique}. But in contrast to NLPs, where SMFCQ is well known as a characterization of the uniqueness of Lagrange multipliers, it is not the case for nonpolyhedral conic programs (including SOCPs), where SRCQ fails to be a characterization of this property; cf.\ \cite[Propositions~4.47 and 4.50]{bs}. As proved in Theorem~\ref{unique}, SRCQ characterizes the uniqueness of Lagrange multipliers for the second-order cone constraint system \eqref{CS} {\em along with} the {\em error bound} estimate \eqref{gf01}, which is automatics for polyhedral systems as in NLPs due to the classical Hoffman lemma. It has been achieved in our proof via the {\em dual qualification condition} \eqref{gf02}, which seems to be new in conic programming while happens to be equivalent to SRCQ in the framework under consideration.

{\bf (ii)} It is worth highlighting the result of Theorem~\ref{unique}(ii) showing that the dual qualification condition \eqref{gf02} yields the full rank of $\nabla\Phi(\ox)$ in \eqref{CS} if $\olm\in\bd(-\Q)\setminus\{0\}$. This is {\em not the case} for NLP constraint systems while reflecting the ``fattiness"  of the second-order cone $\Q$.

{\bf (iii)} Note that the equivalence between \eqref{gf02} and \eqref{sqc} holds true if we replace $\Q$ with any closed {\em convex} sets that is {\em ${\cal C}^2$-cone reducible} in the sense of \cite[Definition~1.135]{bs}. This can be shown by observing that the left-hand side of \eqref{sqc} is convex in this case, and therefore it has a nonempty relative interior in finite dimensions: cf.\ the proof of \cite[Proposition~2.97]{bs}. Note also that Theorem~\ref{unique} can be extended to any {\em second-order regular} convex set $Q$ in the sense of \cite[Definition~3.85]{bs} with the corresponding modifications of the error bound estimate \eqref{gf01}. It is beyond the scope of this paper to provide a proof for such a general framework, and thus we postpone it to our future publications.}
\end{Remark}\vspace*{-0.1in}

To proceed further, define the mapping $\H\colon\R^n\times\R^{m+1}\to\R^{n\times n}$ by
\begin{equation}\label{eqq1}
\H(x;\lambda):=\begin{cases}
-\dfrac{\lambda_0}{\Phi_0(x)}\nabla\widehat{\Phi}(x)^*\nabla\Phi(x)&\textrm{if }\;\Phi(x)=\big(\Phi_0(x),\Phi_r(x)\big)\in\bd(\Q)\setminus\{0\},\\
0&\textrm{otherwise},
\end{cases}
\end{equation}
where $x\in\Gamma$, $\lm=(\lm_0,\lm_r)\in\R\times\R^m$, and $\nabla\widehat{\Phi}(x)=(-\nabla\Phi_0(x),\nabla\Phi_r(x))$. This form is a simplification of the one used in \cite{BoR05,mor}, reflects a nonzero {\em curvature} of the second-order cone $\Q$ at boundary points, and thus is not needed for polyhedra. Recall that $\nabla\Phi(x)$ is an $(m+1)\times n$ matrix and hence $\nabla\widehat{\Phi}(x)^*\nabla\Phi(x)$ is an $n\times n$ matrix in \eqref{eqq1}.

In our derivation of the formula for calculating the graphical derivative of the normal cone mapping $N_\Gamma$ in Section~5, we appeal to the {\em linear conic optimization problem}
\begin{equation}\label{LP}
\min_{\lambda\in\R^{m+1}}\big\{-\left\langle v,\big(\nabla^2\langle\lambda,\Phi\rangle(\bar x)+\H(x;\lambda)\big)v\right\rangle\big|\;\nabla\Phi(\bar x)^*\lambda=\bar x^*\,\;\;\mbox{and}\;\;
\lambda\in N_\Q\big(\Phi(\bar x)\big)\big\}
\end{equation}
generated by the second-order cone $\Q$, where $(\bar x,\bar x^*)\in\gph N_\Gamma$ and $v\in K(\bar x,\bar x^*)$. Denote by $\Lambda(\ox,\ox^*;v)$ the set of optimal solutions to \eqref{LP}. The following result shows that if the primal problem \eqref{LP} has an optimal solution, then  its dual problem has an {\em approximate} feasible solution for which the optimal values of the primal and dual problems are ``almost the same." This is one of the {\em principal differences} between the polyhedral case with the exact duality therein and the nonpolyhedral ice-cream setting. The duality result obtained below is known in case (LMS1) of Proposition~\ref{lag10} (actually in this setting we have  the exact duality; see, e.g., \cite[Theorem~4.14]{ru}), but even in this case our proof is new.\vspace*{-0.1in}

\begin{Theorem}{\bf(approximate duality in linear second-order cone optimization).}\label{Prop1} Taking $(\bar x,\bar x^*)\in\gph N_\Gamma$ and $v\in K(\bar x,\bar x^*)$, suppose that $ \Lambda(\ox,\ox^*)\ne\emp$ and $\Phi(\ox)=0$. Then for every $\tilde{\lm}\in\Lambda(\ox,\ox^*;v)$ and any small $\varepsilon>0$ there exists $z_\varepsilon\in\mathbb{R}^n$ for which we have the relationships
\begin{equation}\label{dua1}
\dist\left(\nabla\Phi(\bar x)z_\varepsilon+\langle v,\nabla^2\Phi(\bar x)v\rangle;\Q\right)\le\varepsilon\quad{and}\quad\;\langle\bar x^*,z_\varepsilon\rangle+\big\langle v,\nabla^2\langle\tilde{\lm},\Phi\rangle(\bar x)v\big\rangle\ge-\varepsilon.
\end{equation}
\end{Theorem}\vspace*{-0.1in}
{\bf Proof.} It follows from \eqref{eqq1} that under $\Phi(\ox)=0$ the optimization problem \eqref{LP} reduces to
\begin{equation}\label{LPv}
\min_{\lambda\in\R^{m+1}}\big\{-\left\langle v,\nabla^2\langle\lambda,\Phi\rangle(\bar x)v\right\rangle
\big|\;\nabla\Phi(\bar x)^*\lambda=\bar x^*\;\;\mbox{and}\;\;\lambda\in-\Q\big\}.
\end{equation}
The {\em dual problem} of \eqref{LPv} can be calculated via \cite[page~125]{bs} and \cite[Example~11.41]{rw} as
\begin{equation}\label{DPv}
\max_{z\in\R^n}\big\{\langle\bar x^*,z\rangle\big|\;\nabla\Phi(\bar x)z+\left\langle v,\nabla^2\Phi(\bar x)v\right\rangle\in T_\Q\big(\Phi(\bar x)\big)\big\}.
\end{equation}
Employing Proposition~\ref{lag10}, we examine all the three possible cases for Lagrange multipliers $\lm\in\Lambda(\ox,\ox^*)$. Picking any $v\in K(\bar x,\bar x^*)$ and $\ve>0$ sufficiently small, consider first the Slater case (LMS1) in Proposition~\ref{lag10} and use the error bound estimate \eqref{error}. This estimate allows us to use the intersection rule from \cite[Proposition~3.2]{io} for the normal cone to $\Lambda(\ox,\ox^*)$  and thus to deduce for any $\tilde{\lm}\in\Lambda(\ox,\ox^*;v)$ that
$$
0\in-\la v,\nabla^2\Phi(\bar x)v\ra+N_{\Lambda(\ox,\ox^*)}(\tilde{\lm})\subset-\la v,\nabla^2\Phi(\bar x)v\ra+N_{-\Q}(\tilde{\lm})+\rge\nabla\Phi(\ox).
$$
This allows us to find some $z\in\R^n$ for which we get
$$
\nabla\Phi(\ox)z+\la v,\nabla^2\Phi(\bar x)v\ra\in N_{-\Q}(\tilde{\lm})\subset\Q=T_\Q\big(\Phi(\ox)\big).
$$
Since $-\Q$ is a convex cone, this inclusion leads us to $\big\la\tilde{\lm},\nabla\Phi(\ox)z+\la v,\nabla^2\Phi(\bar x)v\big\ra=0$. Hence
$$
\la\ox^*,z\ra=\la\tilde{\lm},\nabla\Phi(\ox)z\ra=-\big\la v,\nabla^2\la\tilde{\lm},\Phi\ra(\bar x)v\big\ra,
$$
which in turns implies that $z$ is an optimal solution for the dual problem \eqref{DPv} and that the optimal values of the primal and dual problems agree. Letting $z_\ve:=z$ justifies the validity of both relationships in \eqref{dua1} in the Slater case (LMS1).

In case (LMS2) of Proposition~\ref{lag10}, the set of Lagrange multipliers is a singleton and so is bounded. Using \cite[Proposition~11.39]{rw} tells us that the optimal values of the primal problem \eqref{LPv} and the dual problem \eqref{DPv} agree. Therefore we arrive at
\begin{equation*}
\sup_{z\in\R^n}\left\{\langle\bar x^*,z\rangle\big|\;\nabla\Phi(\bar x)z+\left\langle v,\nabla^2\Phi(\bar x)v\right\rangle\in T_\Q\big(\Phi(\ox)\big)\right\}
=-\left\langle v,\nabla^2\langle\lambda,\Phi\rangle(\bar x)v\right\rangle
\end{equation*}
that allows us for any $\varepsilon>0$ to find $z_\varepsilon$ satisfying the second condition in \eqref{dua1} together with
$$
\nabla\Phi(\bar x)z_\ve+\left\langle v,\nabla^2\Phi(\bar x)v\right\rangle\in T_\Q\big(\Phi(\ox)\big)=\Q.
$$
Thus $z_\varepsilon$ satisfies the first condition in \eqref{dua1} as well, which completes the proof in case (LMS2).

Consider finally case (LMS3) in Proposition~\ref{lag10} where there is $\olm\in\bd(-\Q)$ such that
\begin{equation*}
\Lambda(\ox,\ox^*)=\ker\nabla\Phi(\bar x)^*\cap(-\Q)=\big\{t\bar\lambda\,\big|\,t\ge 0\big\}.
\end{equation*}
In this case the primal problem \eqref{LPv} can be equivalently written as
\begin{equation}\label{dua3}
\min_{\lm\in\R^{m+1}}\left\{-\left\langle v,\nabla^2\langle\lambda,\Phi\rangle(\bar x)v\right\rangle\big|\;\lambda=\alpha\bar\lambda,\,\alpha\ge 0\right\}.
\end{equation}
Since $\Lambda(\ox,\ox^*;v)\ne\emp$, we arrive at $\left\langle v,\,\nabla^2\langle\bar\lambda,\,\Phi\rangle(\bar x)v\right\rangle\leq 0$. Examine the two possible situations:

{\bf(1)} $\left\langle v,\nabla^2\langle\bar\lambda,\Phi\rangle(\bar x)v\right\rangle<0$. In this setting problem \eqref{dua3} has a unique optimal solution $\lambda=0$. Using the arguments similar to the case (LMS2) and applying again \cite[Proposition~11.39]{rw}, we can find some $z_\varepsilon$ satisfying both relationships in \eqref{dua1}.

{\bf(2)} $\left\langle v,\nabla^2\langle\bar\lambda,\Phi\rangle(\bar x)v\right\rangle=0$. In this setting the set of optimal solutions to problem \eqref{dua3} is the entire ray $\{t\bar\lambda\,|\, t\ge 0\}$. Consider now a modified version of (\ref{LPv}) defined by
\begin{equation}\label{LPv2}
\min_{\lambda=(\lm_0,\lm_r)\in\R\times\R^m}\left\{-\left\langle v,\nabla^2\langle\lambda,\Phi\rangle(\bar x)v\right\rangle\big|\;\nabla\Phi(\bar x)^*\lambda=0,\;\lambda\in-Q,\;{-\lm_0\le 1}\right\}.
\end{equation}
Since $\lm\in-\Q$, we get $\|\lm_r\|\le-\lm_0$. This implies that the feasible region of problem (\ref{LPv2}) is nonempty and bounded, and so is
the set of its optimal solutions. Moreover,  its optimal value is zero due to $\left\langle v,\nabla^2\langle\bar\lambda,\Phi\rangle(\bar x)v\right\rangle=0$.
It follows from \cite[Theorem~11.39(a)]{rw} that the optimal value of the dual problem of (\ref{LPv2}) given by
\begin{equation}\label{DPv2}
\max_{(z,\al)\in\R^n\times\R}\left\{\langle 0,z\rangle-\alpha\big|\;\nabla\Phi(\bar x)z+(\alpha,0,\ldots,0)+\left\langle v,\nabla^2\Phi(\bar x)v\right\rangle\in\Q,\;\alpha\ge 0\right\}
\end{equation}
is zero as well. Thus we arrive at the equality
\begin{equation*}
\sup_{(z,\alpha)\in\R^n\times\R}\left\{-\alpha\,\big|\,\nabla\Phi(\bar x)z+(\alpha,0,\ldots,0)+\left\langle v,\nabla^2\Phi(\bar x)v\right\rangle\in\Q,\;\alpha\ge 0\right\}=0.
\end{equation*}
This tells us that for any $\varepsilon>0$ there exists a feasible solution $(z_\varepsilon,\alpha_\varepsilon)\in\mathbb{R}^n\times\mathbb{R}$ to \eqref{DPv2} such that
$-\alpha_\varepsilon>-\ve$. Therefore we have the estimates
\begin{equation*}
\dist\left(\nabla\Phi(\bar x)z_\varepsilon+\langle v,\nabla^2\Phi(\bar x)v\rangle;\Q\right)\le\Vert(\alpha_\varepsilon,0,\ldots,0)\Vert=\alpha_\varepsilon<\varepsilon,
\end{equation*}
which verify the first condition in \eqref{dua1}. Since $\bar x^*=0$ and $\left\langle v,\nabla^2\langle\olm,\Phi\rangle(\bar x)v\right\rangle=0$, we get the second condition in \eqref{dua1} and thus complete the proof of the theorem. $\h$\vspace*{0.05in}

We conclude this section by deriving a second-order sufficient condition for strict local minima in SOCPs needed in what follows. Consider the problem
\begin{equation}\label{socp}
\min\;\ph_0(x)\quad\quad\mbox{subject to}\quad x\in\Gamma,
\end{equation}
where $\ph_0\colon\R^n\to\R$ is twice differentiable, and where $\Gamma$ is taken from \eqref{CS}. Such a second-order sufficient condition was established in \cite[Theorem~3.86]{bs} under the validity of the Robinson constraint qualification \eqref{gf07} that is  equivalent to the metric regularity of the mapping $x\mapsto\Phi(x)-\Q$. It occurs that the same result holds under weaker assumptions on the latter mapping including the validity of MSCQ that guarantees the existence of Lagrange multipliers. \vspace*{-0.1in}

\begin{Proposition}{\bf(second-order sufficient condition for strict local minimizers in SOCP).}\label{sosc} Let $\ox\in\Gamma$ be a feasible solution to \eqref{socp} with $\Phi(\ox)=0$, and let MSCQ hold at $\ox$ while ensuring that $\Lambda(\ox,\ox^*)\ne\emp$ for $\ox^*:=-\nabla\ph_0(\ox)$. Taking any $\olm\in\Lambda(\ox,\ox^*)$, impose the so-called second-order sufficient condition $($SOSC$)$ for optimality:
\begin{equation}\label{sosc1}
\la\nabla^2_{xx}L(\ox,\olm)u,u\ra>0\quad\quad\mbox{for all}\;\;0\ne u\in\big\{u\in\R^n\big|\;\nabla\Phi(\ox)u\in\Q\cap\{\olm\}^\bot\big\},
\end{equation}
where $L(x,\lambda):=\ph_0(x)+\la\lambda,\Phi(x)\ra$. Then $\ox$ is indeed a strict local minimizer for problem \eqref{socp}.
\end{Proposition}\vspace*{-0.1in}
{\bf Proof}. It follows from \cite[Theorem~4.1]{hjo} that MSCQ is as a constraint qualification in \eqref{socp}, i.e., ensures the existence of Lagrange multipliers. Suppose now by that $\ox$ is not a strict local minimizer for \eqref{socp} and thus find a sequence $x_k\to\ox$ as $k\to\infty$ with $\Phi(x_k)\in\Q$ and $\ph_0(x_k)<\ph_0(\ox)$; hence $x_k\ne\ox$. Define $u_k:=\frac{x_k-\ox}{\|x_k-\ox\|}$ and assume without loss of generality that $u_k\to\ou$ for some $0\ne\ou\in\R^n$. It tells us that
$$
\nabla\Phi(\ox)\ou\in\Q\quad\mbox{and}\quad\la\nabla\ph_0(\ox),\ou\ra\le 0.
$$
Combining this with $\olm\in\Lambda(\ox,-\nabla\ph_0(\ox))$ yields $\nabla\Phi(\ox)\ou\in\Q\cap \{\olm\}^\bot$. It is not hard to see that
$$
\ph_0(x_k)-\ph_0(\ox)+\la\olm,\Phi(x_k)\ra\le 0,
$$
which implies by the twice differentiability of $\ph_0$ and $\Phi$ at $\ox$ that
$$
\la\nabla^2_{xx}L(\ox,\olm)\ou,\ou\ra\le 0\quad\quad \mbox{with}\;\;\;\ou\ne 0.
$$
This contradicts \eqref{sosc1} and hence completes the proof of the proposition.$\h$\vspace*{-0.1in}

\section{Graphical Derivative of the Normal Cone Mapping}\sce\vspace*{-0.05in}

Here we present the main result of the paper on calculating the graphical derivative of the normal cone mapping generated by the constraint system \eqref{CS} under imposing {\em merely} the {\em MSCQ} condition. Great progress in this direction was recently made by Gfrerer and Outrata \cite{go16} (preprint of 2014) who calculated this second-order object for polyhedral/NLP constraint systems under MSCQ and a certain additional condition instead of the standard nondegeneracy and Mangasarian-Fromovitz constraint qualifications. Then the additional condition to MSCQ was relaxed in \cite{gm15} and fully dropped subsequently by Chieu and Hien \cite{ch} in the NLP setting. Various calculating formulas for the graphical derivative of the normal cone mappings to nonpolyhedral (including ice-cream) constraints were derived in \cite{go17,mor15,mor}. However, all these results were obtained under some nondegeneracy (a conic extension of the classical linear independence of constraint gradients in NLPs). Thus the graphical derivative formula for the second-order cone constraints given in the next theorem is new even under the Robinson constraint qualification. Furthermore, our proof of this result is significantly different in the major part from that in \cite{go16} and the subsequent developments for polyhedral systems; see Remark~\ref{diss-gd} for more discussions.\vspace*{-0.1in}

\begin{Theorem}{\bf(graphical derivative of the normal cone mapping for the second-order cone constraint systems).}\label{tanf} Let $(\bar x,\bar x^*)\in\gph N_\Gamma$, and let MSCQ from Definition~{\rm\ref{defmscq}} hold at $\ox$ with modulus $\kappa$. Then the tangent cone to $\gph N_\Gamma$ is represented by
\begin{equation}\label{TtN}
\begin{array}{ll}
T_{\small{\gph N_\Gamma}}(\bar x,\bar x^*)=\big\{(v,v^*)\,\big|&v^*\in\big(\nabla^2\langle\lambda,\Phi\rangle(\bar x)+\H(\bar x;\lambda)\big)v+N_{K(\bar x,\bar x^*)}(v)\\&\mbox{for some }\;\lm\in\Lambda(\bar x,\bar x^*;v)\big\},
\end{array}
\end{equation}
where $\Lambda(\bar x,\bar x^*;v)$ is the set of optimal solutions to \eqref{LP} with $\H$ defined in \eqref{eqq1}. Consequently, for all $v\in\R^n$ we have the graphical derivative formula
\begin{equation}\label{gdcs}
(DN_\Gamma)(\bar x,\bar x^*)(v)=\big\{\big(\nabla^2\langle\lambda,\Phi\rangle(\bar x)+\H(\bar x;\lambda)\big)v\,\big|\;\lm\in\Lambda(\bar x,\bar x^*;v)\big\}+N_{K(\bar x,\bar x^*)}(v).
\end{equation}
\end{Theorem}\vspace*{-0.1in}
{\bf Proof}. It is sufficient to justify the tangent cone formula \eqref{TtN}, which immediately yields the graphical derivative one \eqref{gdcs} by definition \eqref{gder}. We split the proof of \eqref{TtN} into three different cases depending on the position of $\Phi(\ox)$ in $\Q$. First assume that $\Phi(\bar x)\in\int(\Q)$ and thus get
\begin{equation*}
\ox^*\in N_\Gamma(\bar x)=\nabla\Phi(\bar x)^*N_\Q(\Phi(\bar x))=\{0\},\quad\quad T_\Gamma(\bar x)=\mathbb{R}^n,\quad\mbox{and}\quad K(\bar x,\bar x^*)=\mathbb{R}^n.
\end{equation*}
By the continuity of $\Phi$ around $\ox$ we find a neighborhood $U$ of $\bar x$ such that $\Phi(x)\in\int(\Q)$ and $N_\Gamma(x)=\{0\}$ whenever $x\in U$.  This tells us that
\begin{equation*}
\gph N_\Gamma\cap[U\times\mathbb{R}^n]=U\times\{0\},
\end{equation*}
which obviously provides the tangent cone representation
\begin{equation}\label{op01}
T_{\small{\gph N_\Gamma}}(\bar x,0)=\mathbb{R}^n\times\{0\}.
\end{equation}
On the other hand, it follows from $\Lambda(\bar x,\bar x^*)=\{0\}$ that $\Lambda(\bar x,\bar x^*;v)=\{0\}$ for all $v\in K(\bar x,\bar x^*)$. This shows that the right-hand side of \eqref{TtN} amounts to $\mathbb{R}^n\times\{0\}$. Combining it with \eqref{op01} verifies the tangent cone formula \eqref{TtN} in this case.

Next we consider the case where $\Phi(\bar x)\in\bd(\Q)\setminus\{0\}$. As argued above, $\Gamma$ can be described in this case by \eqref{CS1} via the mapping $\psi$ from \eqref{mappsi}. Using Lemma~\ref{mscq2} confirms that the mapping $x\mapsto\psi\circ\Phi(x)$ is metrically subregular at $\ox$. Thus it follows from \cite[Theorem~1]{go16} that
\begin{equation}\label{eqq3}
T_{\small{\gph N_\Gamma}}(\bar x,\bar x^*)=\big\{(v,v^*)\,\big|\,v^*\in\nabla^2\langle\tilde\lambda,\,\psi\circ\Phi\rangle(\bar x)v+N_{K(\bar x,\bar x^*)}(v) \;\mbox{for some}\;\tilde\lambda\in\tilde\Lambda(\bar x,\bar x^*;v)\big\},
\end{equation}
where $\tilde\Lambda(\bar x,\bar x^*;v)$ is the set of optimal solutions to the linear program
\begin{equation*}\label{eqq2}
\min_{\tilde\lambda\in\R^2}\big\{-\big\langle v,\nabla^2\langle\tilde\lambda,\psi\circ\Phi\rangle(\bar x)v\big\rangle\,\big|\,\nabla(\psi\circ\Phi)(\bar x)^*\tilde\lambda=\bar x^*,\;\tilde\lambda\in N_{\mathbb{R}^2_-}\big(\psi\circ\Phi(\bar x)\big)\big\}.
\end{equation*}
Define the set of Lagrange multipliers for the  modified constraint system \eqref{CS1} by
$$
\tilde\Lambda(\ox,\ox^*)=\big\{\tilde\lambda\in\R_-^2\,\big|\,\nabla(\psi\circ\Phi)(\bar x)^*\tilde\lambda=\bar x^*,\;\tilde\lambda\in N_{\mathbb{R}^2_-}\big(\psi\circ\Phi(\bar x)\big)\big\}.
$$
It is not hard to observe the implication
\begin{equation}\label{op02}
\tilde\lambda\in\tilde\Lambda(\ox,\ox^*)\Longrightarrow\lambda:=\nabla\psi\big(\Phi(\bar x)\big)^*\tilde\lambda\in\Lambda(\ox,\ox^*),
\end{equation}
where $\Lambda(\ox,\ox^*)$ is taken from \eqref{lagn}. Conversely, we claim that
\begin{equation}\label{op03}
\lambda=(\lm_0,\lm_r)\in\Lambda(\ox,\ox^*)\Longrightarrow\tilde\lambda:=\Big(-\dfrac{\lambda_0}{2\Phi_0(\bar x)},0\Big)\in\tilde\Lambda(\ox,\ox^*).
\end{equation}
To verify \eqref{op03}, we need to show that any $\lambda=(\lm_0,\lm_r)\in\Lambda(\ox,\ox^*)$ can be represented as $\lambda=\nabla\psi(\Phi(\bar x))^*\tilde\lambda$ with some $\tilde\lambda\in N_{\mathbb{R}^2_-}(\psi\circ\Phi(\bar x))$. Since $\Phi(\bar x)=(\Phi_0(\ox),\Phi_r(\ox))\in\bd(\Q)\setminus\{0\}$, it follows that $(\psi\circ\Phi)(\bar x)=(0,-\Phi_0(\bar x))$ and $\Phi_0(\bar x)>0$, which lead us to $N_{\mathbb{R}^2_-}((\psi\circ\Phi)(\bar x))=\mathbb{R}_+\times\{0\}$. Thus we  need to find some $\alpha\ge 0$ such that the pair $\tilde\lambda=(\alpha,0)$ satisfies the equation
\begin{equation*}
\lambda=\nabla\psi\big(\Phi(\bar x)\big)^*\tilde\lambda=\begin{bmatrix}
-2\Phi_0(\bar x)&-1\\2\Phi_r(\bar x)&0\\
\end{bmatrix}\begin{pmatrix}
\alpha\\0\\
\end{pmatrix}=2\alpha\widehat{\Phi}(\bar x),
\end{equation*}
which is clearly fulfilled for $\alpha=-\dfrac{\lambda_0}{2\Phi_0(\bar x)}$ and hence justifies the claimed implication \eqref{op03}.
Using these observations brings us to the following relationships:
\begin{eqnarray*}
\nabla^2\langle\tilde\lambda,\psi\circ\Phi\rangle(\bar x)&=&\nabla^2\left(\alpha(-\Phi_0^2(\cdot)+\Vert\Phi_r(\cdot)\Vert^2)\right)(\bar x)=2\alpha\nabla\left[\widehat{\Phi}(\cdot)^*\nabla\Phi(\cdot)\right](\bar x)\\
&=&2\alpha\left[\nabla\widehat{\Phi}(\bar x)^*\nabla\Phi(\bar x)+\langle\widehat{\Phi}(\bar x),\nabla^2\Phi(\bar x)\rangle\right]\\
&=&2\alpha\nabla\widehat{\Phi}(\bar x)^*\nabla\Phi(\bar x)+\langle 2\alpha\widehat{\Phi}(\bar x),\nabla^2\Phi(\bar x)\rangle\\
&=&-\dfrac{\lambda_0}{\Phi_0(\bar x)}\nabla\widehat{\Phi}(\bar x)^*\nabla\Phi(\bar x)+\langle\lambda,\nabla^2\Phi(\bar x)\rangle=\H(\bar x;\lambda)+\nabla^2\langle \lambda,\Phi\rangle(\bar x).
\end{eqnarray*}
Combining it with \eqref{op02} and \eqref{op03} confirms that \eqref{eqq3} reduces to \eqref{TtN} in this case.

It remains to consider the most difficult nonpolyhedral case where $\Phi(\ox)=0$. We begin with verifying the inclusion ``$\subset$" in \eqref{TtN}. Picking any $(v,v^*)\in T_{\small{\gph N_\Gamma}}(\bar x,\bar x^*)$, observe that it suffices to show the validity of the following two inclusions:
\begin{equation}\label{rhs}
v\in K(\bar x,\bar x^*)\quad\mbox{and}\quad v^*-\nabla^2\langle\bar\lambda,\Phi\rangle(\bar x)v\in N_{K(\bar x,\bar x^*)}(v)\;\;\;
\mbox{for some}\;\;\bar\lambda\in\Lambda(\bar x,\bar x^*;v).
\end{equation}
To proceed, we get from the tangent cone definition \eqref{2.5} that for $(v,v^*)\in T_{\small{\gph N_\Gamma}}(\bar x,\bar x^*)$ there are sequences
$t_k\downarrow 0$ and $(v_k,v_k^*)\to(v, v^*)$ as $k\to\infty$ such that
$$
(x_k,x_k^*):=(\bar x+t_kv_k,\bar x^*+t_kv_k^*)\in\gph N_\Gamma,\quad k\in\N.
$$
Let us split the subsequent proof of the inclusion ``$\subset$" in \eqref{TtN} into the four steps.

\textbf{Step~1:} {\em There exists a sequence $\{\lambda_k\in\Lambda(x_k,x_k^*)\}$ with $\lambda_k\to\bar\lambda$ as $k\to\infty$ for some $\bar\lambda\in\Lambda(\bar x,\bar x^*)$.} To verify this statement, we deduce first directly from \cite[Lemma~2.1]{gm16} and the robustness of MSCQ that there is a positive number $\delta$ such that $x_k\in\Gamma\cap\B_\dd(\ox)$ and that
\begin{equation*}
\Lambda(x_k,x_k^*)\cap\kappa\Vert x_k^*\Vert\B\ne\emp\quad\mbox{for all}\;\;k\in\N,
\end{equation*}
where $\kappa>0$ is the constant taken from Definition~{\rm\ref{defmscq}}. This allows us to find $\lambda_k\in\Lambda(x_k,x_k^*)$ so that $\Vert\lambda_k\Vert\le\kappa\Vert x_k^*\Vert$ for all $k\in\N$. Thus the boundedness of $\{x_k^*\}$ yields the one for $\{\lambda_k\}$, and therefore $\lambda_k\to\bar\lambda$ for some $\olm\in\Rm$ along a subsequence.  In this way we conclude that $\bar\lambda\in\Lambda(\bar x,\bar x^*)$, where the latter set is represented by \eqref{lag0} due to $\Phi(\ox)=0$.

\textbf{Step~2:} {\em We have $v\in T_\Gamma(\bar x)\cap\{\bar x^*\}^\perp=K(\bar x,\bar x^*)$}. The equality here is by the definition of the critical cone \eqref{critco}; so getting the first one in \eqref{rhs} requires only the verification the claimed inclusion. To furnish this, recall first from \eqref{first2} that $T_\Gamma(\bar x)=\left\{w\in\mathbb{R}^{n}\,\big|\,\nabla\Phi(\bar x)w\in\Q\right\}$. It follows from $x_k\in\Gamma$ for all $k\in\N$ and $\Phi(\ox)=0$ that
\begin{equation*}
\Phi(x_k)=t_k\nabla\Phi(\bar x)v_k+o(t_k)\in\Q,\quad k\in\N.
\end{equation*}
Dividing the latter by $t_k$ and passing to the limit as $k\to\infty$ yield $\nabla\Phi(\bar x)v\in\Q$, and so $v\in T_\Gamma(\bar x)$. Since $\bar\lambda\in\Lambda(\bar x,\bar x^*)$ and $\langle\lambda_k,\Phi(x_k)\rangle=0$ for all $k\in\N$, we get
\begin{equation*}
\begin{array}{lll}
\langle\bar x^*,v\rangle&=&\langle\nabla\Phi(\bar x)^*\bar\lambda,v\rangle=\langle\bar\lambda,\nabla\Phi(\bar x)v\rangle=\disp\lim_{k\to\infty}\langle\lambda_k,\nabla\Phi(\bar x)v_k\rangle\\&=&\disp\lim_{k\to\infty}\dfrac{\langle\lambda_k,\Phi(x_k)+o(t_k)\rangle}{t_k}=\lim_{k\to\infty}\Big\langle\lambda_k,\dfrac{o(t_k)}{t_k}\Big\rangle=0
\end{array}
\end{equation*}
and thus finish the proof of the statement in Step~2.

\textbf{Step~3:} {\em We have the inclusion $v^*-\nabla^2\langle\bar\lambda,\Phi\rangle(\bar x)v\in\big(K(\bar x,\bar x^*)\big)^*$ for the multiplier $\bar\lambda\in\Lambda(\bar x,\bar x^*)$ constructed in Step~{\rm 1}.} Indeed, by the definition of $x_k^*$ we get
\begin{equation*}
v_k^*=\dfrac{x_k^*-\bar x^*}{t_k}=\dfrac{\nabla\Phi(x_k)^*\lambda_k-\bar x^*}{t_k}=\dfrac{\nabla\Phi(\bar x)^*\lambda_k+t_k\nabla^2\langle\lambda_k,\Phi\rangle(\bar x)v_k+o(t_k)-\bar x^*}{t_k},
\end{equation*}
which in turn leads us to the equality
\begin{equation}\label{eq10}
v_k^*-\nabla^2\langle\lambda_k,\Phi\rangle(\bar x)v_k+\dfrac{o(t_k)}{t_k}
=\nabla\Phi(\bar x)^*\dfrac{\lambda_k}{t_k}-\dfrac{\bar x^*}{t_k}.
\end{equation}
Using $\lambda_k\in-\Q=N_\Q(\Phi(\bar x))$ and \eqref{eq8} yields $v_k^*-\nabla^2\langle\lambda_k,\Phi\rangle(\bar x)v_k+\dfrac{o(t_k)}{t_k}\in\left(K(\bar x,\bar x^*)\right)^*$.
Since $\left(K(\bar x,\bar x^*)\right)^*$ is closed, the passage to to the limit as $k\to\infty$ gives us the desired inclusion.

\textbf{Step~4:} {\em We have $\bar\lambda\in\Lambda(\bar x,\bar x^*;v)$ and $\left\langle v,v^*-\nabla^2\langle\bar\lambda,\Phi\rangle(\bar x)v\right\rangle=0$ for the multiplier $\bar\lambda$ constructed above.} To furnish this, it suffices to show that
\begin{equation}\label{eq9}
\left\langle v,\nabla^2\langle\lambda,\Phi\rangle(\bar x)v\right\rangle\le\left\langle v,\nabla^2\langle\bar\lambda,\Phi\rangle(\bar x)v\right\rangle\;\mbox{ for any }\;\lambda\in\Lambda(\bar x,\bar x^*).
\end{equation}
Picking $\lambda\in\Lambda(\bar x,\bar x^*)$ gives us $\lambda\in-\Q$ by \eqref{lag0}. Using this together with $\Phi(x_k)\in\Q$ and $\langle\lambda_k,\Phi(x_k)\rangle=0$, we get the relationships
\begin{eqnarray*}
0&\le&-\langle\lambda,\Phi(x_k)\rangle=\langle\lambda_k-\lambda,\Phi(x_k)\rangle\\
&=&t_k\langle\lambda_k-\lambda,\nabla\Phi(\bar x)v_k\rangle+\dfrac{1}{2}t_k^2\left\langle v_k,\nabla^2\langle\lambda_k-\lambda,\Phi\rangle(\bar x)v_k\right\rangle+o(t_k^2)\\
&=&t_k\langle\nabla\Phi(\bar x)^*\lambda_k-\bar x^*,v_k\rangle+\dfrac{1}{2}t_k^2\left\langle v_k,\nabla^2\langle\lambda_k-\lambda,\Phi\rangle(\bar x)v_k\right\rangle+o(t_k^2).
\end{eqnarray*}
Dividing by $t_k^2$ and employing \eqref{eq10} bring us to
\begin{equation*}
0\le\big\langle v_k,v_k^*-\nabla^2\langle\lambda_k,\Phi\rangle(\bar x)v_k+\dfrac{o(t_k)}{t_k}\big\rangle+\dfrac{1}{2}\left\langle v_k,\nabla^2\langle\lambda_k-\lambda,\Phi\rangle(\bar x)v_k\right\rangle+\dfrac{o(t_k^2)}{t_k^2},
\end{equation*}
which implies by passing to the limit as $k\to\infty$ that
\begin{equation}\label{dua10}
0\le\left\langle v,v^*-\nabla^2\langle\bar\lambda,\Phi\rangle(\bar x)v\right\rangle+\dfrac{1}{2}\left\langle v,\nabla^2\langle\bar\lambda-\lambda,\Phi\rangle(\bar x)v\right\rangle.
\end{equation}
It follows from the relationships proved in Steps~2 and 3 that
\begin{equation}\label{eq7}
\left\langle v,v^*-\nabla^2\langle\bar\lambda,\Phi\rangle(\bar x)v\right\rangle\le 0.
\end{equation}
which together with \eqref{dua10} yields \eqref{eq9}. Finally, since \eqref{dua10} holds for any $\lambda\in\Lambda(\bar x,\bar x^*)$, letting $\lambda=\bar\lambda$ therein results in the inequality
\begin{equation*}
\left\langle v,v^*-\nabla^2\langle\bar\lambda,\Phi\rangle(\bar x)v\right\rangle\ge 0.
\end{equation*}
Combining it with \eqref{eq7} justifies Step~4, and thus we arrive at the inclusion ``$\subset$" in \eqref{TtN}.

Now we give a detailed proof of the opposite inclusion in \eqref{TtN}, which occurs to be more involved. Pick $(v,v^*)$ from the right-hand side of \eqref{TtN}, which satisfies \eqref{rhs} in the case of $\Phi(\ox)=0$ under consideration. We proceed by showing that there are sequences $t_k\downarrow 0$ and $x_k\to\bar x$ as $k\to\infty$ satisfying the conditions
\begin{equation}\label{eq4}
\bar x+t_kv-x_k=o(t_k)\quad\mbox{and}\quad\dist\big(\bar x^*+t_kv^*;N_\Gamma(x_k)\big)=o(t_k),\quad k\in\N.
\end{equation}
These guarantee the existence of $x_k^*\in N_\Gamma(x_k)$ such that
$$
(x_k,x_k^*)=\Big(\ox+t_k\big(v+\frac{o(t_k)}{t_k}\big),\ox^*+t_k\big(v^*+\frac{o(t_k)}{t_k}\big)\Big)\in\gph N_\Gamma,
$$
and thus we arrive at $(v,v^*)\in T_{\small{\gph N_\Gamma}}(\bar x,\bar x^*)$, which is the goal.

To furnish it, we conclude by the choice of $(v,v^*)$ and the usage of Theorem~\ref{norco} that there are $\lambda\in\Lambda(\bar x,\bar x^*;v)$ and $\mu\in T_{-\Q}(\lambda)$ satisfying the equalities
\begin{equation}\label{dua9}
v^*=\nabla^2\langle\lambda,\Phi\rangle(\bar x)v+\nabla\Phi(\bar x)^*\mu\quad\mbox{and}\quad\langle\mu,\nabla\Phi(\bar x)v\rangle=0.
\end{equation}
It comes from $\mu\in T_{-\Q}(\lambda)$ that there are sequences $t_i\downarrow 0$ and $\mu_i\to\mu$ as $i\to\infty$ with $\lambda+t_i\mu_i\in-\Q$.
Choose $\alpha>0$ so small that $\alpha\Vert\nabla^2\langle\lambda,\Phi\rangle(\bar x)\Vert\le\dfrac{1}{2}$ holds. This ensures that the matrix $I+\alpha\nabla^2\langle\lambda,\Phi\rangle(\bar x)$ is positive-definite, where $I$ is the $n\times n$ identity matrix. Proposition~\ref{sosc} tells us that there exists $r>0$ such that $\bar x$ is the strict global minimizer for the problem
\begin{equation}\label{eq2}
\min_{x\in\R^n}\left\{\Vert\bar x+\alpha\bar x^*-x\Vert^2\,\big|\,x\in\Gamma\cap\B_r(\ox)\right\}.
\end{equation}
For any fixed $k\in\N$ we select a positive number $\varepsilon_k<\big(16\alpha k^2(\kappa\Vert\bar x^*\Vert+1)\big)^{-1}$. Since $\lambda$ solves the linear optimization problem \eqref{LPv}, Theorem~\ref{Prop1} ensures the existence of $z_k\in\mathbb{R}^n$ with
\begin{equation}\label{dua6}
{\rm dist}\left(\nabla\Phi(\bar x)z_k+\left\langle v,\nabla^2\Phi(\bar x)v\right\rangle;\Q\right)\le\varepsilon_k\quad\mbox{and}\quad\langle\bar x^*,z_k\rangle
+\left\langle v,\nabla^2\langle\lambda,\Phi\rangle(\bar x)v\right\rangle\ge-\varepsilon_k.
\end{equation}
Picking next $i\in\N$, consider yet another optimization problem
\begin{equation}\label{eq3}
\min_{x\in\R^n}\big\{\Vert\bar x+t_iv+\dfrac{1}{2}t_i^2z_k +\alpha(\bar x^*+t_iv^*)-x\Vert^2\,\big|\,x\in\Gamma\cap\B_r(\ox)\big\},
\end{equation}
which admits an optimal solution due to the classical Weierstrass theorem. This global minimizer $x_i$ is surely unique for each $i$, and it is not hard to check that $x_i\to\bar x$ as $i\to\infty$. Indeed, suppose that $x_i\to\tilde x$ for some $\tilde x$ along a subsequence, we see that
\begin{equation*}
\Vert\bar x+\alpha\bar x^*-\tilde x\Vert^2\le\Vert\bar x+\alpha\bar x^*-x\Vert^2\quad\mbox{for all}\;\;x\in\Gamma\cap\B_r(\ox),
\end{equation*}
which yields $\tilde x=\bar x$ since $\ox$ is the strict global minimizer for \eqref{eq2}. Assume now without loss of generality that $x_i\in\int\B_r(\ox)$ for all $i\in\N$ and utilize the first-order necessary optimality condition from \cite[Proposition~5.1]{m06} at $x_i$ for problem \eqref{eq3} to get the following inclusion:
\begin{equation}\label{dua8}
\alpha(\bar x^*+t_iv^*)+t_i\Big(\dfrac{\bar x+t_iv-x_i}{t_i}+\dfrac{1}{2}t_iz_k\Big)\in N_\Gamma(x_i).
\end{equation}
It follows from $\Phi(\ox)=0$ and the twice differentiability of $\Phi$ around $\ox$ that
\begin{equation*}
\Phi(\bar x+t_iv+\dfrac{1}{2}t_i^2z_k)=t_i\nabla\Phi(\bar x)v+\dfrac{1}{2}t_i^2\Big((\nabla\Phi(\bar x)z_k+\left\langle v,\nabla^2\Phi(\bar x)v\right\rangle\Big)+o(t_{i}^2).
\end{equation*}
Since $v$ satisfies \eqref{rhs}, we get $\nabla\Phi(\bar x)v\in T_\Q(\Phi(\ox))=\Q$. Taking this into account along with the first inequality in \eqref{dua6}, the latter equality yields the estimate
\begin{equation*}
\dist\Big(\Phi(\bar x+t_iv+\dfrac{1}{2}t_i^2z_k);\Q\Big)\le\dfrac{\varepsilon_k}{2}t_i^2+o(t_i^2),
\end{equation*}
which together with the assumed MSCQ at $\ox$ results in
\begin{equation*}
\dist\Big(\bar x+t_iv+\dfrac{1}{2}t_i^2z_k;\Gamma\Big)\le\dfrac{\kappa\varepsilon_k}{2}t_i^2+o(t_i^2).
\end{equation*}
This guarantees that for any $i\in\N$ there exists $\tilde x_i\in\Gamma$ such that
\begin{equation}\label{eq5}
\Vert\bar x+t_iv+\dfrac{1}{2}t_i^2z_k-\tilde x_i\Vert\le\dfrac{\kappa\varepsilon_k}{2}t_i^2+o(t_i^2),
\end{equation}
and so we verify that $\tilde x_i\to\bar x$ as $i\to\infty$. This tells us that $\tilde x_i\in\Gamma\cap\B_r(\ox)$ for all $i$ sufficiently large. Since $x_i$ is a global minimizer for \eqref{eq3}, we get
\begin{equation*}
\big\Vert\bar x+t_iv+\dfrac{1}{2}t_i^2z_k+\alpha(\bar x^*+t_iv^*)-x_i\big\Vert^2\le\big\Vert\bar x+t_iv+\dfrac{1}{2}t_i^2z+\alpha(\bar x^*+t_iv^*)-\tilde x_i\big\Vert^2
\end{equation*}
for all large $i$, which together with \eqref{eq5} leads us to the estimates
\begin{eqnarray*}
&&\big\Vert\bar x+t_iv+\dfrac{1}{2}t_i^2z_k-x_i\big\Vert^2+2\alpha\big\langle\bar x^*+t_iv^*,\bar x+t_iv+\dfrac{1}{2}t_i^2z_k-x_i\big\rangle\\
&\le&\big\Vert\bar x+t_iv+\dfrac{1}{2}t_i^2z_k-\tilde x_i\big\Vert^2+2\alpha\big\langle\bar x^*+t_iv^*,\bar x+t_iv+\dfrac{1}{2}t_i^2z_k-\tilde x_i\big\rangle\\
&\le&\big\Vert\bar x+t_iv+\dfrac{1}{2}t_i^2z_k-\tilde x_i\big\Vert^2+2\alpha\left(\Vert\bar x^*\Vert+t_i\Vert v^*\Vert\right)\big\Vert\bar x+t_iv+\dfrac{1}{2}t_i^2z_k-\tilde x_i\big\Vert\\
&\le&\alpha\kappa\Vert\bar x^*\Vert\varepsilon_k t_i^2+o(t_i^2).
\end{eqnarray*}
These yield in turn the relationships
\begin{eqnarray}\label{eq5.1}
\big\Vert\bar x+t_iv+\dfrac{1}{2}t_i^2z_k-x_i\big\Vert^2&\le&-2\alpha\big\langle\bar x^*+t_iv^*,\bar x+t_iv+\dfrac{1}{2}t_i^2z_k-x_i\big\rangle+\alpha\kappa\Vert \bar x^*\Vert\varepsilon_k t_i^2+o(t_i^2)\nonumber\\
&=&2\alpha\big[\langle\bar x^*+t_iv^*,x_i-\bar x\rangle-t_i\langle\bar x^*,v\rangle-t_i^2\langle v^*,v\rangle-\dfrac{1}{2}t_i^2\langle\bar x^*,z_k\rangle\big]\nonumber\\
&&+\alpha\kappa\Vert\bar x^*\Vert\varepsilon_k t_i^2+o(t_i^2).
\end{eqnarray}
Recall further from the first inclusion in \eqref{rhs} that $v\in K(\bar x,\bar x^*)$ and hence $\langle\bar x^*,v\rangle=0$. It follows from \eqref{dua9} and \eqref{dua6}, respectively, that
\begin{equation}\label{eq5.2}
\langle v^*,v\rangle=\left\langle v,\nabla^2\langle\lambda,\Phi\rangle(\bar x)v\right\rangle\quad\mbox{and}\quad-\langle\bar x^*,z_k\rangle\le\left\langle v,\nabla^2\langle\lambda,\Phi\rangle(\bar x)v\right\rangle+\varepsilon_k.
\end{equation}

Next we are going to find an upper estimate for the first term on the right-hand side of the equality in \eqref{eq5.1}. It follows from both equalities in \eqref{dua9} that
\begin{eqnarray*}
\langle\bar x^*+t_iv^*,x_i-\bar x\rangle&=&\left\langle\nabla\Phi(\bar x)^*\lambda+t_i\left(\nabla^2\langle\lambda,\Phi\rangle(\bar x)v+\nabla\Phi(\bar x)^*\mu\right),x_i-\bar x\right\rangle\\&=&\left\langle\lambda+t_i\mu,\nabla\Phi(\bar x)(x_i-\bar x)\right\rangle+t_i\left\langle v,\nabla^2\langle\lambda,\Phi\rangle(\bar x)(x_i-\bar x)\right\rangle\\&=&\left\langle\lambda+t_i\mu_i,\nabla\Phi(\bar x)(x_i-\bar x)\right\rangle+t_i\left\langle\mu-\mu_i,\nabla\Phi(\bar x)(x_i-\bar x)\right\rangle\nonumber\\&&+t_i\left\langle v,\nabla^2\langle\lambda,\Phi\rangle(\bar x)(x_i-\bar x)\right\rangle\\&=&\big\langle\lambda+t_i\mu_i,\Phi(x_i)-\dfrac{1}{2}\langle x_i-\bar x,\nabla^2\Phi(\bar x)(x_i-\bar x)\rangle\big\rangle\\
&&+\;t_i\left\langle v,\nabla^2\langle\lambda,\Phi\rangle(\bar x)(x_i-\bar x)\right\rangle+o(t_i\Vert x_i-\bar x\Vert)+o(\Vert x_i-\bar x\Vert^2)\nonumber\\&=&\langle\lambda+t_i\mu_i,\Phi(x_i)\rangle-\dfrac{1}{2}\left\langle x_i-\bar x,\nabla^2\langle\lambda,\Phi\rangle(\bar x)(x_i-\bar x)\right\rangle\\&&+t_i\left\langle v,\nabla^2\langle\lambda,\Phi\rangle(\bar x)(x_i-\bar x)\right\rangle+o(t_i\Vert x_i-\bar x\Vert)+o(\Vert x_i-\bar x\Vert^2).
\end{eqnarray*}
Using these together with $\lambda+t_i\mu_i\in-\Q$ and $\Phi(x_i)\in\Q$ brings us to the estimate
\begin{eqnarray}\label{eq5.4}
\langle\bar x^*+t_iv^*,x_i-\bar x\rangle&\le&-\dfrac{1}{2}\left\langle\bar x-x_i,\nabla^2\langle\lambda,\Phi\rangle(\bar x)(\bar x-x_i)\right \rangle\nonumber-t_i\big\langle v,\nabla^2\langle\lambda,\Phi\rangle(\bar x)(\bar x-x_i)\big\rangle\\&&+o(t_i\Vert x_i-\bar x\Vert )+o(\Vert x_i-\bar x\Vert^2).
\end{eqnarray}
Combining now the conditions in \eqref{eq5.1}--\eqref{eq5.4}, we arrive at the following relationships:
\begin{eqnarray*}
\big\Vert\bar x+t_iv+\dfrac{1}{2}t_i^2z_k-x_i\big\Vert^2&\le&2\alpha\big[-\dfrac{1}{2}\left\langle\bar x-x_i,\nabla^2\langle\lambda,\Phi\rangle(\bar x)(\bar x-x_i)\right\rangle-t_i\left\langle v,\nabla^2\langle\lambda,\Phi\rangle(\bar x)(\bar x-x_i)\right\rangle\big.\\
&&\hspace{0.26in}-t_i^2\left\langle v,\nabla^2\langle\lambda,\Phi\rangle(\bar x)v\right\rangle+\dfrac{1}{2}t_i^2\left(\left\langle v,\nabla^2\langle\lambda,\Phi\rangle(\bar x)v\right\rangle+\varepsilon_k\right)\big]\\&&\hspace{0.26in}+o(t_i\Vert x_i-\bar x\Vert)+o(\Vert x_i-\bar x\Vert^2)+\alpha\kappa\Vert\bar x^*\Vert\varepsilon_k t_i^2+o(t_i^2)\\&=&-\alpha\big[\left\langle\bar x-x_i,\nabla^2\langle\lambda,\Phi\rangle(\bar x)(\bar x-x_i)\right\rangle+2t_i\left\langle v,\nabla^2\langle\lambda,\Phi\rangle(\bar x)(\bar x-x_i)\right\rangle\\
&&\hspace{0.26in}+t_i^2\left\langle v,\nabla^2\langle\lambda,\Phi\rangle(\bar x)v\right\rangle\big]+\alpha\varepsilon_kt_i^2+\alpha\kappa\Vert\bar x^*\Vert\varepsilon_k t_i^2\\&&\hspace{0.26in}+o(t_i\Vert x_i-\bar x\Vert)+o(\Vert x_i-\bar x\Vert^2)+o(t_i^2)\\
&=&-\alpha\left\langle\bar x+t_iv-x_i,\nabla^2\langle\lambda,\Phi\rangle(\bar x)(\bar x+t_iv-x_i)\right\rangle+\alpha(\kappa\Vert\bar x^*\Vert+1)\varepsilon_kt_i^2\\&&\hspace{0.26in}+o(t_i\Vert x_i-\bar x\Vert)+o(\Vert x_i-\bar x\Vert^2)+o(t_i^2)\\
&\le&\dfrac{1}{2}\Vert\bar x+t_iv-x_i\Vert^2+\alpha(\kappa\Vert\bar x^*\Vert+1)\varepsilon_kt_i^2\\
&&\hspace{0.26in}+o(t_i\Vert x_i-\bar x\Vert)+o(\Vert x_i-\bar x\Vert^2)+o(t_i^2),
\end{eqnarray*}
where the last inequality comes from the fact that the matrix $\dfrac{1}{2}I+\alpha\nabla^2\langle\lambda,\Phi\rangle(\bar x)$ is positive-semidefinite.
This allows us to conclude that
\begin{equation*}
\big\Vert\bar x+t_iv+\dfrac{1}{2}t_i^2z_k-x_i\big\Vert^2-\dfrac{1}{2}\Vert\bar x+t_iv-x_i\Vert^2\le\alpha(\kappa\Vert\bar x^*\Vert+1)\varepsilon_kt_i^2+o(t_i\Vert x_i-\bar x\Vert)+o(\Vert x_i-\bar x\Vert ^2)+o(t_i^2),
\end{equation*}
which verifies the validity of the inequality
\begin{eqnarray*}
\dfrac{1}{2}\Vert\bar x+t_iv-x_i\Vert^2+t_i^2\langle z_k,\bar x-x_i\rangle+t_i^3\langle z_k,v\rangle+\dfrac{1}{4}t_i^4\Vert z_k\Vert^2&\le&\alpha(\kappa\Vert\bar x^*\Vert+1)\varepsilon_kt_i^2+o(t_i\Vert x_i-\bar x\Vert)\\
&&+o(\Vert x_i-\bar x\Vert^2)+o(t_i^2).
\end{eqnarray*}
Since $\varepsilon_k<\dfrac{1}{16\alpha(\kappa\Vert\bar x^*\Vert+1)}$, the latter inequality can be simplified as
\begin{eqnarray}\label{eq6}
\Vert\bar x+t_iv-x_i\Vert^2&\le&2\alpha(\kappa\Vert\bar x^*\Vert+1)\varepsilon_kt_i^2+o(t_i\Vert x_i-\bar x\Vert)+o(\Vert x_i-\bar x\Vert ^2)+o(t_i^2)\\
&\le&\dfrac{1}{8}t_i^2+o(t_i^2)+o(t_i\Vert x_i-\bar x\Vert)+o(\Vert x_i-\bar x\Vert^2),\nonumber
\end{eqnarray}
and therefore we get for all $i$ sufficiently large that
\begin{equation*}
\Vert\bar x+t_iv-x_i\Vert\le\dfrac{1}{2}\left(t_i+\Vert x_i-\bar x\Vert\right).
\end{equation*}
In this way we arrive at the estimates
\begin{equation*}
\Vert x_i-\bar x\Vert\le\Vert\bar x+t_iv-x_i\Vert+t_i\Vert v\Vert\le\dfrac{1}{2}t_i+\dfrac{1}{2}\Vert x_i-\bar x\Vert,
\end{equation*}
which in turn imply that $\Vert x_i-\bar x\Vert=O(t_i)$ and so $o(t_i\Vert x_i-\bar x\Vert)=o(\Vert x_i-\bar x\Vert^2)=o(t_i^2)$.
Using these relationships together with \eqref{eq6} gives us
\begin{equation*}
\Vert\bar x+t_iv-x_i\Vert^2\le 2\alpha(\kappa\Vert\bar x^*\Vert+1)\varepsilon_kt_i^2+o(t_i^2),
\end{equation*}
and so we come by passing to the limit as $i\to\infty$ to the inequalities
\begin{equation*}
\lim_{i\to\infty}\dfrac{\Vert\bar x+t_iv-x_i\Vert^2}{t_i^2}\le 2\alpha(\kappa\Vert\bar x^*\Vert+1)\varepsilon_k\le\frac{1}{8k^2}.
\end{equation*}
Remember that $k\in\N$ has been fixed through the above proof of the inclusion ``$\supset$" in \eqref{TtN}. This allows us to find an index $i_k$ for which we have the estimates
\begin{equation}\label{final}
\dfrac{\Vert\bar x+t_{i_k}v-x_{i_k}\Vert}{t_{i_k}}\le\dfrac{1}{2k}\quad\quad\mbox{and}\quad\quad t_{i_k}\Vert z_k\|\le\dfrac{1}{k}.
\end{equation}
Repeating this process for any $k\in\N$, we construct sequences $t_{i_k}$ and $x_{i_k}$ that satisfy \eqref{final} and such that
$t_{i_k}\downarrow 0$ and $x_{i_k}\to\ox$ as $k\to\infty$.
Combining finally \eqref{final} and \eqref{dua8} leads us to
\begin{equation*}
\dfrac{\dist\big(\bar x^*+t_{i_k}v^*;N_\Gamma(x_{i_k})\big)}{t_{i_k}}\le\dfrac{1}{k}.
\end{equation*}
It yields \eqref{eq4} with $t_k:=t_{i_k}$ and $x_k:=x_{i_k}$ and so completes the proof of the theorem.$\h$\vspace*{0.03in}

It is worth mentioning an equivalent neighborhood version of the pointbased formula \eqref{TtN} in Theorem~\ref{tanf}, which is an ice-cream counterpart of the polyhedral result established recently by Gfrerer and Ye \cite[Theorem~4]{gy}.\vspace*{-0.1in}

\begin{Corollary}{\bf(neighborhood representation of the tangent cone to the normal cone graph for ice-cream constraint systems).}\label{fuzzy}
Under the assumptions of Theorem~{\rm\ref{tanf}} there is $\dd>0$ such that for all $x\in\Gamma\cap \B_\dd(\ox)$ and all $x^*\in N_\Gamma(x)$ we have the representations
\begin{equation}\label{TtN2}
\begin{array}{llll}
T_{\small{\gph N_\Gamma}}(x,x^*)&=&\big\{(v,v^*)\,\big|&v^*\in\big(\nabla^2\langle\lambda,\Phi\rangle(x)+\H(x;\lambda)\big)v+N_{K( x,x^*)}(v)\\
&& &\mbox{for some}\;\;\lm\in\Lambda(x,x^*;v)\cap\kappa\|x^*\|\B\big\},
\end{array}
\end{equation}
\begin{equation}\label{gp}
(DN_\Gamma)(x, x^*)(v)=\big\{\big(\nabla^2\langle\lambda,\Phi\rangle(x)+\H(x;\lambda)\big)v\big|\;\lm\in\Lambda(x,x^*;v)\cap\kappa\|x^*\|\B\big\}+N_{K(x,x^*)}(v).
\end{equation}
\end{Corollary}\vspace*{-0.1in}
{\bf Proof.} It is clear that \eqref{TtN2} yields \eqref{TtN}. To verify the reverse implication, pick $\kappa,\dd>0$ from Step~1 in the proof of Theorem~\ref{tanf} and suppose by the robustness of MSCQ that it holds at any $x\in\Gamma\cap\B_\dd(\ox)$. If $\Lambda(x,x^*;v)=\emp$, then both sides in \eqref{TtN2} are empty. Otherwise, we proceed as Steps~1-4 in the proof of Theorem~\ref{TtN} to establish the following inclusion:
\begin{equation*}
\begin{array}{llll}
T_{\small{\gph N_\Gamma}}(x,x^*)&\subset&\big\{(v,v^*)\,\big|&v^*\in \big(\nabla^2\langle\lambda,\Phi\rangle(x)+\H( x;\lambda)\big)v+N_{K(x,x^*)}(v)\\
&& &\mbox{for some}\;\;\lm\in\Lambda(x,x^*;v)\cap\kappa\|x^*\|\B\big\}.
\end{array}
\end{equation*}
On the other hand, it is proved in Theorem~\ref{TtN} that the set
\begin{equation*}
\big\{(v,v^*)\,\big|\,v^*\in\big(\nabla^2\langle\lambda,\Phi\rangle( x)+\H(x;\lambda)\big)v+N_{K(x,x^*)}(v)\;\;\mbox{for some}\;\;\lm\in\Lambda(x,x^*;v)\big\}\\
\end{equation*}
is contained in $T_{\small{\gph N_\Gamma}}(x,x^*)$. Having all of this, we arrive at the claimed equivalence. The obtained representation \eqref{TtN2} yields the graphical derivative one \eqref{gp} by its definition. $\h$\vspace*{0.05in}

The next consequence of Theorem~\ref{tanf} concerns an important case of the tangent cone formula in the case where $\ox^*=0$, which is used in what follows.\vspace*{-0.1in}

\begin{Corollary} {\bf(simplification of the graphical derivative formula for $\ox^*=0$).}\label{tanf2} Let $\ox^*=0$ in the framework of Theorem~{\rm\ref{tanf}}. Then we have
\begin{equation}\label{eq25}
T_{\small{\gph N_\Gamma}}(\bar x,0)=\big\{(v,v^*)\,\big|\,v^*\in N_{K(\ox,0)}(v)\big\}=\gph N_{K(\ox,0)}
\end{equation}
and correspondingly the graphical derivative formula
\begin{equation}\label{eq26}
(DN_\Gamma)(\bar x,0)(v)=N_{K(\bar x,0)}(v)=\nabla\Phi(\bar x)^*\big[N_\Q\big(\Phi(\ox)\big)\cap\{\nabla\Phi(\bar x)v\}^\bot\big].
\end{equation}
\end{Corollary}\vspace*{-0.1in}
{\bf Proof}. If $\ox^*=0$, we deduce from \eqref{TtN2} that $\lm=0$. Using this together with $\H(\bar x;\lm)=0$ for $\lm=0$, we arrive at \eqref{eq25} and hence at \eqref{eq26}.  $\h$\vspace*{-0.1in}

\begin{Remark}{\bf (extensions to products of second-order cones).}\label{product}
{\rm It is not hard to observe that Theorem~\ref{tanf} can be extended to the finite products of the second-order cones. Indeed, consider the extended second-order cone constraint system defined by
$$
\Gamma:=\big\{x\in\R^n\,\big|\,\Phi(x):=\big(\Phi^1(x),\ldots,\Phi^J(x)\big)\in\Q\big\},
$$
where $\Phi^{j}\colon\R^n\to\R^{{m_j}+1}$ as $j=1,\ldots,J$ are twice differentiable, and where the set $\Q$ is given by
\begin{equation} \label{cone Q}
\Q:=\prod_{j=1}^{J}{\cal Q}_{m_j+1}
\end{equation}
via the second-order cones $\Q_{m_j+1}$ in $\R^{m_j+1}$. Then it is easy to see that
$$
\big((v^1,\ldots,v^J),(w^1,\ldots,w^J)\big)\in\gph N_\Gamma\Longleftrightarrow(v^j,w^j)\in\gph N_{\Gamma_{m_j+1}}\quad\mbox{for}\;\;j=1,\ldots,J,
$$
where $\Gamma_{m_j+1}:=\big\{x\in\R^n\,\big|\,\Phi^{m_j+1}(x)\in{\Q_{m_j+1}}\big\}$. Employing now \cite[Proposition~6.41]{rw} allows us to derive the desired counterpart of \eqref{TtN} for the product constraint system $\Gamma$.}
\end{Remark}\vspace*{-0.2in}

\begin{Remark}{\bf(discussions on the graphical derivative formulas).}\label{diss-gd}{\rm

{\bf (i)} First we highlight some important differences between our proof of Theorem~\ref{tanf} for {\em nonpolyhedral} second-order constraint systems and its {\em polyhedral} counterpart for NLPs in \cite[Theorem~1]{go16} and in the similar devices from \cite{ch,gm15}. Unlike the latter proof that exploits the {\em Hoffman lemma} to verify the inclusion ``$\subset$" in \eqref{TtN}, we {\em do not appeal} to any error bound estimate; this is new even for polyhedral systems. Our approach is applicable to other cone-constrained frameworks; however, we believe that some {\em error bound} estimate {\em is needed} for the general setting. The reason for avoiding error bounds in the proof of Theorem~\ref{tanf} is that in the ice-cream case we have the inclusion $N_\Q(x)\subset N_\Q(0)$ for any $x\in\Rm$. Another difference between our proof and that in \cite{go16} lies in the justification of the inclusion ``$\supset$" in the tangent cone formula. Indeed, the proof in \cite{go16} employs the {\em exact duality}, which holds in the polyhedral setting. In contrast, our proof relies on the {\em approximate duality} established in Theorem~\ref{Prop1}.

{\bf (ii)} The first result on the tangent cone and the graphical derivative of normal cone mapping to the general conic constraint system
\begin{equation}\label{conc}
\Gamma:=\big\{x\in\R^n\big|\;\Phi(x)\in\Th\big\},
\end{equation}
where $\Th\subset\R^m$ is a closed and convex, was established by Mordukhovich, Outrata ana Ram\'irez \cite[Theorem~3.3]{mor15} under the {\em nondegeneracy condition} from \cite{bs} and the rather restrictive assumption on the convexity of $\Gamma$. This result was derived not in the form of \eqref{TtN} but in terms of the directional derivative of the {\em projection mapping} associated with $\Th$. Later the same authors improved this result in \cite[Theorem~5.2]{mor} by dropping the convexity of $\Gamma$ under the {\em projection derivation condition} discussed in Section~\ref{sec3}, which enabled them to write the main result for \eqref{conc} in the form of \eqref{TtN}. However, as proved in Corollary~\ref{pdpn}, this PDC does not hold at nonzero boundary points of $\Q$ and so \cite[Theorem~5.2]{mor}---obtained also under the nondegeneracy condition ---cannot be utilized in the ice-cream framework when $\Phi(\ox)\in\bd(\Q)\setminus\{0\}$.

{\bf (iii)} Quite recently, Gfrerer and Outrata \cite[Theorem~2]{go17} calculated the graphical derivative of the normal mapping to \eqref{conc} under the validity of the {\em nondegeneracy condition} when $\Th$ is not necessarily convex. Combining their result with Corollary~\ref{dcod} above in the ice-cream framework, we see that it agrees with Theorem~\ref{TtN} {\em provided that} the nondegeneracy condition is satisfied. However, our results can be applied to much broader settings since it only demands the fulfillment of MSCQ. As mentioned above, our results seem to be new for SOCPs even under the validity of RCQ \eqref{gf07}, which is equivalent to the metric regularity of $x\mapsto\Phi(x)-\Q$ around $(\ox,0)$. Note that in the latter case the Lagrange multiplier set $\Lambda(\ox,\ox^*)$ admits either the (LMS1) or the (LMS2) representation from its description in Proposition~\ref{lag0}.}
\end{Remark}\vspace*{-0.2in}

\section{Examples and Application to Isolated Calmness}\label{sec05}\sce\vspace*{-0.05in}

First we illustrate the applicability of the main result in Theorem~\ref{tanf} to the ice-cream constraint systems at points where neither nondegeneracy nor Robinson constraint qualification is satisfied.\vspace*{-0.1in}

\begin{Example}{\bf(calculation of graphical derivative for ice-cream normal cone systems).}\label{Ex1} {\rm
Define the mapping $\Phi:\mathbb{R}^2\to\mathbb{R}^3$ by
\begin{equation*}
\Phi(x):=\big(\sqrt{2}x_1^2+x_2,x_1^2+\dfrac{1}{\sqrt{2}}x_2,x_1^2-\dfrac{1}{\sqrt{2}}x_2\big)\quad\mbox{for}\quad x=(x_1,x_2)\in\R^2
\end{equation*}
and consider the constraint system associated with the three-dimensional ice-cream cone $\Q_3$:
\begin{equation*}
\Gamma=\big\{x\in\mathbb{R}^2\,\big|\,\Phi(\bar x)\in\Q_3\big\}=\big\{(x_1,x_2)\in\mathbb{R}^2\,\big|\,x_2\ge 0\big\}.
\end{equation*}
Given any $x\in\Gamma$, we claim that the mapping $x\mapsto\Phi(x)-\Q_3$ is metrically subregular at $(x,0)$, i.e., MSCQ holds at $x$. To begin with, observe by \eqref{bg01} and direct calculations that
\begin{equation*}
\dist\big((x_1,x_2);\Gamma \big)=\begin{cases}
0\quad&\textrm{if }\ x_2\ge 0,\\
-x_2&\textrm{if }\;x_2<0;
\end{cases}
\end{equation*}
\begin{equation*}
\dist\big( \Phi(x_1,x_2);\Q_3\big)=\begin{cases}
0\quad&\textrm{if }\; x_2\ge 0,\\
-\sqrt{2}x_2&\textrm{if }\;x_1=0,\,x_2<0,\\
\dfrac{\sqrt{2}}{2}\Big(-x_2+\sqrt{2x_1^4+x_2^2}-\sqrt{2}x_1^2\Big)\quad&\textrm{otherwise},
\end{cases}
\end{equation*}
which gives us $\dist((x_1,x_2);\Gamma)\le\sqrt{2}\dist(\Phi(x_1, x_2);\Q_3)$ for all $(x_1,x_2)\in\mathbb{R}^2$ and thus verifies the validity of
MSCQ at any $x\in \Gamma$. It is not hard to check that
\begin{equation*}
N_\Gamma(x)=\left\{\begin{array}{lll}
\{(0,0)\}&\textrm{if }\;x_2>0,\\
\{0\}\times\mathbb{R}_-\quad&\textrm{if }\;x_2=0,\\
\varnothing&\textrm{if }\;x_2<0
\end{array}\right.
\quad\mbox{and}\quad
T_\Gamma(x)=\left\{\begin{array}{lll}
\mathbb{R}^2&\textrm{if }\;x_2>0,\\
\mathbb{R}\times\mathbb{R}_+\quad&\textrm{if }\;x_2=0,\\
\varnothing&\textrm{if }\;x_2<0.
\end{array}
\right.
\end{equation*}
On the other hand, the direct calculation tells us that
\begin{equation}\label{gh03}
T_{\small{\gph N_\Gamma}}(\bar x,\bar x^*)=\begin{cases}
[\mathbb{R}\times(0,\infty)\times\{(0,0)\}]\cup[\mathbb{R}\times\{0\}\times\{0\}\times\mathbb{R}_-]&\textrm{if}\;x_2= 0,\bar x^*=0,\\
\mathbb{R}\times\{0\}\times\{0\}\times\mathbb{R}
&\textrm{if}\;x_2= 0,\bar x^*\ne 0,\\
\mathbb{R}^2\times\{(0,0)\}
&\textrm{if}\;x_2>0,\bar x^*=0.
\end{cases}
\end{equation}

Let us now apply Theorem~\ref{tanf} to calculate the tangent cone to $\gph N_\Gamma$ and the graphical derivative of the normal cone mapping. For $\lm=(\lm_0,\lm_1,\lm_2)\in\R^3$ we have
\begin{eqnarray*}
\nabla\Phi(x)^*&=&\begin{bmatrix}
2\sqrt{2}x_1&2x_1&2x_1\\
1&\dfrac{1}{\sqrt{2}}&-\dfrac{1}{\sqrt{2}}\\
\end{bmatrix},\quad\quad
\nabla^2\langle\lambda,\Phi\rangle(x)=\begin{bmatrix}
2\sqrt{2}\lambda_0+2\lambda_1+2\lambda_2&0\\
0&0\\
\end{bmatrix}.
\end{eqnarray*}
Consider further the following five characteristic cases:

{\bf Case~1:} $\bar x=(0, 0)$ and $\bar x^*=(0,0)\in N_\Gamma(\bar x)$. In this case we have $\Phi(\bar x)=0$, $\H(\bar x;\lm)=0$, and $K(\bar x,\bar x^*)=T_\Gamma(\ox)=\R\times\R_+$.  Applying Corollary~\ref{tanf2} tells us that
$$
T_{\small{\gph N_\Gamma}}(\bar x,\bar x^*)=\gph N_{K(\bar x,\bar x^*)}=\big[\mathbb{R}\times(0,\infty)\times\{(0,0)\}\big]\cup\big[\mathbb{R}\times\{0\}\times\{0\}\times\mathbb{R}_-\big],
$$
$$
(DN_\Gamma)(\bar x,\bar x^*)\big((v_1,v_2)\big)=N_{K(\bar x,\bar x^*)}\big((v_1,v_2)\big)=\begin{cases}
\big\{(0,0)\big\}&\textrm{if }\;v_2>0,\\
\big\{0\big\}\times\R_-&\textrm{if }\;v_2=0
\end{cases}
$$
for $v=(v_1,v_2)$, which agrees with the calculation in \eqref{gh03}.

{\bf Case~2:} $\ox=(0,0)$ and $\bar x^*=(0,-1)$ with $K(\bar x,\bar x^*)=\mathbb{R}\times\{0\}$. Take $\big((v_1,v_2),(v_1^*,v^*_2)\big)$ from the right-hand side of \eqref{TtN} and observe that for any $v:=(v_1,v_2)\in K(\bar x,\bar x^*)$ it holds
$$
N_{K(\bar x,\bar x^*)}(v)=\{0\}\times\mathbb{R}\quad\mbox{and}\quad\Lambda(\ox,\ox^*;v)=
\begin{cases}
\Big\{(-1,\dfrac{1}{\sqrt{2}},\dfrac{1}{\sqrt{2}})\Big\}&\textrm{if }\;v_1\ne 0,\\
\big\{\lambda\in-\Q_3\,\big|\,\sqrt{2}\lambda_0+\lambda_1-\lambda_2=-\sqrt{2}\big\}&\textrm{if }\;v_1=0.
\end{cases}
$$
Thus Theorem~\ref{tanf} gives us the following inclusions:
\renewcommand{\labelenumi}{{\rm(\roman{enumi})}}
\begin{enumerate}
\item if $v_1\ne 0$ and $v_2=0$, then
\begin{equation*}
v^*\in\begin{bmatrix}
-2\sqrt{2}+\sqrt{2}+\sqrt{2}&0\\
0&0\\
\end{bmatrix}\begin{pmatrix}
v_1\\
0\\
\end{pmatrix}+\{0\}\times\mathbb{R}=\big\{0\big\}\times\mathbb{R};
\end{equation*}

\item if $v_1=v_2=0$, then there exists $\lambda\in\Lambda(\ox,\ox^*;v)$ such that
\begin{equation*}
v^*\in\begin{bmatrix}
2\sqrt{2}\lambda_1+2\lambda_2+2\lambda_3&0\\
0&0\\
\end{bmatrix}\begin{pmatrix}
0\\
0\\
\end{pmatrix}+\{0\}\times\mathbb{R}=\big\{0\big\}\times\mathbb{R}.
\end{equation*}
\end{enumerate}
We therefore arrive at the tangent cone formula
\begin{equation*}
T_{\small{\gph N_\Gamma}}(\bar x,\bar x^*)=\big\{(v,v^*)\,\big|\,v_2=0\;\;\;\mbox{and}\;\;\;v_1^*=0\,\big\},
\end{equation*}
which yields for $v=(v_1,v_2)$ with $v_2=0$ the graphical derivative one
$$
(DN_\Gamma)(\bar x,\bar x^*)\big((v_1,v_2)\big)=\big\{0\big\}\times\R.
$$
Thus in this case we again agree with the calculation in \eqref{gh03}.

{\bf Case~3:} $\ox=(1,0)$ and $\ox^*=(0,0)\in N_\Gamma(\ox)$. Observe that in this case we have $\Phi(\bar x)\in\bd(\Q_3)\setminus\{0\}$, $K(\ox,\ox^*)=\R\times\R_+$, and it follows from \eqref{eqq1} that
\begin{equation*}
\H(x;\lambda)=-\dfrac{\lambda_0}{\sqrt{2}}\begin{bmatrix}
0&-2\sqrt{2}\\
-2\sqrt{2}&0\\
\end{bmatrix}.
\end{equation*}
Applying Corollary~\ref{tanf2} gives us the same formulas for $T_{\small{\gph N_\Gamma}}$ and $DN_\Gamma$ as in Case~1.

{\bf Case~4:} $\ox=(1,0)$ and $\bar x^*=(0,-1)\in N_\Gamma(\bar x)$ with $K(\bar x,\bar x^*)=\mathbb{R}\times\{0\}$. Taking $\big((v_1,v_2),(v_1^*,v^*_2)\big)$ from the right-hand side of \eqref{TtN}, observe that for all $v=(v_1,v_2)\in K(\bar x,\bar x^*)$ we get $N_{K(\bar x,\bar x^*)}(v)=\{0\}\times\mathbb{R}$. It is easy to check that
$$
\Lambda(\ox,\ox^*)=\Lambda(\ox,\ox^*;v)=\Big\{\Big(-1,\dfrac{1}{\sqrt{2}},\dfrac{1}{\sqrt{2}}\Big)\Big\},
$$
which implies that for any $\lm\in\Lambda(\ox,\ox^*;v)$ we have
\begin{equation*}
\nabla^2\langle\lambda,\Phi\rangle(\bar x)+\H(\bar x;\lambda)=
\begin{bmatrix}
-2\sqrt{2}+\sqrt{2}+\sqrt{2}&0\\
0&0\\
\end{bmatrix}+\dfrac{1}{\sqrt{2}}\begin{bmatrix}
0&-2\sqrt{2}\\
-2\sqrt{2}&0\\
\end{bmatrix}=\begin{bmatrix}
0&-2\\
-2&0\\
\end{bmatrix}.
\end{equation*}
Appealing to Theorem~\ref{tanf} tells us that
\begin{eqnarray*}
T_{\small{\gph N_\Gamma}}(\bar x,\bar x^*)&=&\big\{\big((v_1,v_2),( v_1^*,v^*_2)\big)\,\big|\,v_2=0,\,v^*\in(0,-2v_1)+\{0\}\times\mathbb{R} \big\}\\
&=&\big\{\big((v_1,v_2),( v_1^*,v^*_2)\big)\,\big|\,v_2=0,\,v_1^*=0\big\},
\end{eqnarray*}
which readily implies that for any $v=(v_1,v_2)$ with $v_2=0$ we get
$$
(DN_\Gamma)(\bar x,\bar x^*)\big((v_1,v_2)\big)=\big\{0\big\}\times\R.
$$

{\bf Case~5:} $\bar x=(0,1)$ and $\bar x^*=(0,0)$. In this case we have $K(\bar x,\bar x^*)=\mathbb{R}^2$ and so $N_{K(\bar x,\bar x^*)}(v)=\{(0,0)\}$ for all $v\in\R^2$. It is easy to see that $\Lambda(\ox,\ox^*)=(\sqrt{2},-1,1)\mathbb{R}_-$, which tells us that the Lagrange multipliers set has the representation in (LMS3) of Proposition~\ref{lag10}. Employing again Corollary~\ref{tanf2} ensures the validity of the relationships
$$
T_{\small{\gph N_\Gamma}}(\bar x,\bar x^*)=\gph N_{K(\bar x,\bar x^*)}=\mathbb{R}^2\times\big\{(0,0)\big\},
$$
and therefore we arrive at the graphical derivative formula
$$
(DN_\Gamma)(\bar x,\bar x^*)((v_1,v_2))=\big\{(0,0)\big\},\quad v\in\R^2,
$$
which illustrates the applicability of Theorem~\ref{tanf} under the imposed MSCQ condition. Since the set of Lagrange multiplies is {\em unbounded} in some cases above, the {\em metric regularity} condition (equivalent to the Robinson constraint qualification, which characterizes the boundedness of Lagrange multipliers) {\em fails} in this example, not even talking about the nondegeneracy condition. This completes our considerations in this example.}
\end{Example}\vspace*{-0.1in}

Next we provide an application of Theorem~\ref{tanf} to an important stability property well recognized in variational analysis and optimization; see, e.g., \cite{dsz,dr,mor15} and the references therein. Recall that a mapping $F\colon\R^n\tto\R^m$ is said to be {\em isolatedly calm} at $(\ox,\oy)\in\gph F$ if there exist a constant $\ell\ge0$ and neighborhoods $U$ of $\ox$ and $V$ of $\oy$ such that
\begin{equation}\label{calm}
F(x)\cap V\subset\{\oy\}+\ell\|x-\ox\|\B\;\mbox{ for all }\;x\in U.
\end{equation}
In what follows we apply the graphical derivative formula established above to characterize the isolated calmness property of the {\em parametric variational system}
\begin{equation}\label{VS}
S(p)=\big\{x\in\R^n\,\big|\,p\in f(x)+N_\Gamma(x)\big\}
\end{equation}
generated by the the ice-cream cone $\Q\subset\Rm$ via \eqref{CS}, where $f\colon\R^n\to\R^n$ is a differentiable mapping. The following theorem provides a complete characterization of the isolated calmness of the variational system \eqref{VS} entirely via its given data.\vspace*{-0.1in}

\begin{Theorem}{\bf(characterization of isolated calmness for ice-cream variational systems).}\label{calm-VS} Let $(\op,\ox)\in\gph S$ with $S$ taken from \eqref{VS}. In addition to the standing assumptions on $\Gamma$ from \eqref{CS} and the MSCQ condition of Theorem~{\rm\ref{tanf}}, suppose that $f$ is Fr\'echet differentiable at $\ox\in\Gamma$. Then $S$ enjoys the isolated calmness property at $(\op,\ox)$ if and only if
\begin{equation}\label{crc5}
\left\{\begin{array}{ll}
0\in\nabla f(\ox)v+\big(\nabla^2\langle\lambda,\Phi\rangle(\bar x)+\H(\bar x;\lambda)\big)v+N_{K(\bar x,\op-f(\ox))}(v)\\
\lm\in\Lambda\big(\ox,\op-f(\ox);v\big)\cap\kappa\,\|\op-f(\ox)\|\B
\end{array}\right.
\Longrightarrow v=0,
\end{equation}
where $\kappa>0$ is the metric subregularity constant of the mapping $x\mapsto \Phi(x)-\Q$ at $(\ox,0)$.
\end{Theorem}\vspace*{-0.1in}
{\bf Proof.} We invoke a graphical derivative characterization of the isolated calmness property \eqref{calm} for arbitrary closed-graph multifunctions written as
\begin{equation}\label{calm-cr}
DF(\ox,\oy)(0)=\{0\}.
\end{equation}
This result goes back to Rockafellar \cite{roc87} although it was not explicitly formulated in \cite{roc87}; see \cite[Theorem~4C.1]{dr} with the commentaries. It easily follows from the Fr\'echet differentiability of $f$ at $\ox$ and the structure of $S$ in \eqref{VS} that $v\in DS(\op,\ox)(u)$ if and only if $u\in\nabla f(\ox)v+(DN_\Gamma)(\ox,\op-f(\ox))(v)$. Using now the calmness criterion \eqref{calm-cr} and substituting there the graphical derivative formula from Corollary~\ref{fuzzy}, we arrive at the claimed characterization \eqref{crc5}.$\h$\vspace*{0.03in}

Finally in this section, we present a numerical example of the ice-cream variational system \eqref{VS} where the application of Theorem~\ref{calm-VS} allows us to reveal that the isolated calmness property holds at some feasible points while failing at other ones.

\begin{Example}{\bf(verification of isolated calmness).}\label{calm-exa}{\rm Consider the variational system \eqref{VS}, where $f\colon\mathbb{R}^2\to\mathbb{R}^2$ is given by
\begin{equation*}
f(x):=\big(x_1,x_2^2\big)\quad\textrm{for}\quad x=(x_1,x_2)
\end{equation*}
while the constraint set $\Gamma$ is taken from Example~\ref{Ex1}.
We examine the following cases:

{\bf Case~1:} $\ox=(0,0)$ and $\bar p=f(\ox)=(0,0)$. In this case we have
\begin{equation*}
\nabla f(\bar x)v+DN_\Gamma\big(\ox,\bar p-f(\ox)\big)\big((v_1,v_2)\big)=\begin{pmatrix}
v_1\\0\\
\end{pmatrix}+\begin{cases}
\{(0,0)\}&\textrm{if }v_2>0,\\
\{0\}\times\R_-\quad&\textrm{if }\;v_2=0.
\end{cases}
\end{equation*}
Invoking the corresponding calculations from Example~\ref{Ex1} shows implication \eqref{crc5} does not hold. Thus the isolated calmness of \eqref{VS} fails at this point $(\op,\ox)$.

{\bf Case~2:} $\ox=(0,0)$ and $\bar p=(0,-1)$. In this case we have $\bar p-f(\ox)=(0,-1)$ and
\begin{equation*}
\nabla f(\bar x)v+DN_\Gamma\big(\ox,\bar p-f(\ox)\big)\big((v_1,v_2)\big)=\begin{pmatrix}
v_1\\0\\
\end{pmatrix}+\{0\}\times\R\quad\textrm{if }\;v_2=0.
\end{equation*}
It is clear that implication \eqref{crc5} holds for this case, and so does the isolated calmness at $(\op,\ox)$.

{\bf Case~3:} $\ox=(1,0)$ and $\bar p=f(\ox)=(1,0)$. The right-hand side of the inclusion in \eqref{crc5} for this case is the same as that in Case~1. Therefore we come up with the same conclusion that isolated calmness does not hold at this point.

{\bf Case~4:} $\ox=(1,0)$ and $\bar p=(1,-1)$. We get the validity of the same implication \eqref{crc5} as that in Case~2 and therefore justify the isolated calmness of \eqref{VS} at the point under consideration.

{\bf Case~5:} $\ox=(0,1)$ and $\bar p=f(\ox)=(0,1)$. Then the right-hand side of the inclusion in \eqref{crc5} reduces to $(v_1,2v_2)+\{(0,0)\}$. It is easy to see that implication \eqref{crc5} holds, which therefore justifies the isolated calmness of \eqref{VS} in this case.}
\end{Example}\vspace*{-0.25in}

\section{Concluding Remarks}\vspace*{-0.05in}

This paper provides a comprehensive second-order analysis of conic constraint systems associated with the second-order/Lorentz/ice-cream cone $\Q$. In particular, it gives precise calculations---entirely via the initial system data---of the graphical derivative of the conic constraint \eqref{CS} when the constraint system is merely metrically subregular. To the best our knowledge, results of this type have been established so far either for polyhedral systems, or under the constraint nondegeneracy that is much stronger than metric subregularity.

In our future research we plan to extend the obtained results to other nonpolyhedral constraint systems in conic programming and to give applications of the established calculations of the graphical derivative to various areas of variational analysis and optimization where this construction and related ones naturally appear.
\end{document}